%% file: main-new.tex
\documentclass[letterpaper, 10 pt, conference]{ieeeconf}  

\IEEEoverridecommandlockouts                              

\usepackage{amsmath} 
\usepackage{amssymb}  
\usepackage{mathtools}
\usepackage{bm,fixmath}
\usepackage{hyperref}
\usepackage{graphicx,tikz}
\usepackage{pgfplots} 
\pgfplotsset{compat=newest}
\pgfplotsset{plot coordinates/math parser=false}
\definecolor{mycolor1}{rgb}{0.00000,0.44700,0.74100}%
\definecolor{mycolor2}{rgb}{0.85000,0.32500,0.09800}%
\definecolor{mycolor3}{rgb}{0.92900,0.69400,0.12500}
\usetikzlibrary{calc,shapes,arrows.meta,decorations.pathreplacing,decorations.markings,positioning}

\usetikzlibrary{positioning}
\tikzset{%
	dot/.style={
   	circle,inner sep=0pt,minimum size=0pt,fill=white
   	},
    block/.style={draw, fill=white, rectangle, 
            minimum height=2em, minimum width=3em},
    blockSmall/.style={draw, fill=white, rectangle, 
            minimum height=1.75em, minimum width=2em},
    input/.style={inner sep=0pt},       
    output/.style={inner sep=0pt},      
    sum/.style = {draw, fill=white, circle, minimum size = 0.1cm, node distance=1.5cm, inner sep=0pt},
    pinstyle/.style = {pin edge={to-,thin,black}}
}
\pgfkeys{/pgf/number format/.cd,1000 sep={\,}}
\newtheorem{assum}{Assumption}
\newtheorem{rem}{Remark}
\DeclareMathOperator*{\argmin}{arg\,min}
\DeclareMathOperator{\tr}{tr}
\newcommand{\norm}[1]{\left\lVert#1\right\rVert}
\newcommand{\abs}[1]{\left\lvert#1\right\rvert}
\def\N{{\mathbb{N}}}

\def\RNumb{{\mathbb{R}}}

\def\C{{\mathbb{C}}}

\def\gradf{{\nabla f_k}}
\def\NDist{{\mathcal{N}}}

\def\zero{{\mathbf{0}}}

\def\c{{\mathbold{c}}}

\def\e{{\mathbold{e}}}

\def\g{{\mathbold{g}}}

\def\w{{\mathbold{w}}}

\def\x{{\mathbold{x}}}

\def\m{{\mathbold{m}}}
\def\A{{\mathbold{A}}}
\def\B{{\mathbold{B}}}
\def\C{{\mathbold{C}}}

\def\F{{\mathbold{F}}}
\def\G{{\mathbold{G}}}
\def\H{{\mathbold{H}}}
\def\I{{\mathbold{I}}}

\def\P{{\mathbold{P}}}
\def\Q{{\mathbold{Q}}}

\def\V{{\mathbold{V}}}

\def\K{{\mathbold{K}}}
\def\J{{\mathbold{J}}}

\def\xiB{{\mathbold{\xi}}}

\def\Fhat{\widehat{\F}}
\def\Fhatext{\widehat{\F}_{ext}}

\def\Ghat{\widehat{\G}}
\def\Ghatext{\widehat{\G}_{ext}}

\def\Hhat{\widehat{\H}}
\def\Hhatext{\widehat{\H}_{ext}}

\def\Jhat{\widehat{J}}

\def\cbarS{\overline{c}}

\def\cbar{\overline{\c}}
\def\cbarS{\overline{c}}

\def\xiBar{\overline{\xiB}}

\def\xiH{\widehat{\xiB}}
\def\Fext{\F_{ext}}
\def\Hext{\H_{ext}}
\def\Gext{\G_{ext}}
\def\Kext{\K_{ext}}
\def\Jext{\J_{ext}}

\def\LB{\mathbold{\Lambda}}

\title{\LARGE \bf
Stochastic models for online optimization
}

\author{Umberto Casti$^{1}$ and Sandro Zampieri$^{2}$
\thanks{*This work was not supported by any organization}
\thanks{$^{1}$Umberto Casti is with the Department of Information Engineering (DEI), University of Padova, 35131 Padova, Italy (\href{mailto:castiumber@dei.unipd.it}{\texttt{castiumber@dei.unipd.it.}})}%
\thanks{$^{2}$Sandro Zampieri is with the Department of Information Engineering (DEI), University of Padova, 35131 Padova, Italy (\href{mailto:zampi@dei.unipd.it}{\texttt{zampi@dei.unipd.it.}})}%
}

\begin{document}

\maketitle
\thispagestyle{empty}
\pagestyle{empty}

\begin{abstract}
    In this paper, we propose control-theoretic methods as tools for the design of online optimization algorithms that are able to address dynamic, noisy, and partially uncertain time-varying quadratic objective functions. Our approach introduces two algorithms specifically tailored for scenarios where the cost function follows a stochastic linear model. The first algorithm is based on a Kalman filter-inspired approach, leveraging state estimation techniques to account for the presence of noise in the evolution of the objective function. The second algorithm applies $\mathcal{H}_\infty$-robust control strategies to enhance performance under uncertainty, particularly in cases in which model parameters are characterized by a high variability.

Through numerical experiments, we demonstrate that our algorithms offer significant performance advantages over the traditional gradient-based method and also over the optimization strategy proposed in~\cite{bastianello_internal_2022} based on deterministic models. 
\end{abstract}


\section{Introduction}\label{sec: introduction}

Online optimization is a rapidly developing field with applications across diverse areas, including control~\cite{liao-mcpherson_semismooth_2018,paternain_realtime_2019}, signal processing~\cite{hall_online_2015,fosson_centralized_2021,natali_online_2021}, and machine learning~\cite{shalev_online_2011,dixit_online_2019,chang_distributed_2020}. These problems have grown in significance due to advances in technology, as they involve the optimization of time-varying cost functions within dynamic environments~\cite{dallanese_optimization_2020,simonetto_timevarying_2020}.

A common approach in online optimization is to adapt static methods, like gradient descent, for their use in a dynamic setting. These approaches, known as \textit{unstructured} methods, generally achieve convergence near the optimal trajectory, but do not leverage the dynamics of the problem to enhance performance~\cite{simonetto_timevarying_2020}.

In contrast, \textit{structured} algorithms aim to exploit the evolution of the cost function to improve tracking performance. By incorporating assumptions about the rate of change of the optimization problem, structured methods can achieve a better accuracy over time. A prominent example of such approaches is the class of so called \textit{prediction-correction} algorithms~\cite{simonetto_dual_2019,bastianello_extrapolation_2023}, which use past information on the cost function to predict future changes, thereby adjusting and correcting the solution as the problem evolves. These methods have demonstrated better tracking performance than \textit{unstructured} approaches by embedding assumptions about the rate of change in the cost function. Specifically, they utilize time-derivative assumptions to predict and adjust the optimizer solution based on past observations. However, these methods do not leverage a more precise model of the cost function's time evolution. In this work, we align our results with the model-based solution proposed in~\cite{bastianello_internal_2022}, which integrates a model for cost evolution without requiring bounded assumptions on the rate of change. This approach employs control-theoretic tools, particularly the Internal Model Principle (IMP)~\cite{fadali_digital_2019}, to design a novel online optimization algorithm for time-varying quadratic cost functions with a linearly evolving term governed by a deterministic LTI system. 

In contrast, this paper extends the methods introduced in~\cite{bastianello_internal_2022} to scenarios in which the time-varying cost function follows a \textit{stochastic} linear model. Within this new framework, we propose two algorithms: the first is inspired by Kalman filtering, and the second employs $\mathcal{H}_\infty$-robust control techniques. We illustrate the effectiveness of these algorithms through numerical experiments.

Before presenting our contributions, we review related work, particularly within the intersection of control theory and optimization. One important class of online optimization problems, related to \textit{feedback optimization}, plays a key role in applications such as Model Predictive Control (MPC). In MPC and similar applications, the output of a dynamic system is fed into an optimization algorithm, which then generates the control input, closing the control loop. Here, optimization techniques are applied to address control problems, whereas in this and other studies, control theory is applied to design and analyze optimization algorithms. For instance, control techniques have been applied in both static and online optimization contexts, as seen in~\cite{lessard_analysis_2016,scherer_optimization_2023,davydov_contracting_2023}. Specifically,~\cite{lessard_analysis_2016,scherer_optimization_2023} focus on static optimization, while~\cite{davydov_contracting_2023} uses contraction theory to analyze continuous-time online optimization algorithms.

\noindent \textbf{Notation.} We denote by $\N$ and $\RNumb$ the sets of natural and real numbers, respectively. Vectors and matrices are denoted by bold letters, \emph{e.g.} $\x\in\mathbb{R}^n$ and $\A \in \mathbb{R}^{n\times n}$. For the $i$-th component of a vector, we use the notation $\left[\x\right]_i$. Conjugate transposition of a vector $\x$ or a matrix $\A$ is denoted by $\x^*$ or $\A^*$, respectively. The identity matrix of dimension $n$ is denoted by $\I_n$, $\zero$ denotes the vectors of all zeros. The $\mathrm{2}$-norm of a vector and the $\mathcal{H}_{2}$-norm of a system are both denoted by $\norm{\cdot}_2$ while the $\mathcal{H}_\infty$-norm of a system is denoted by $\norm{\cdot}_{\infty}$. For matrices symmetric $\A$ and $\B$, the notation $\A \preceq \B$ indicates that $\B - \A$ is positive semi-definite. The symbol $\otimes$ denotes the Kronecker product.
For a function $f: \RNumb^n \to \RNumb$, $\nabla f$ denotes its gradient. Given a matrix $\A$, we denote by $\sigma_{\max} \left(\A\right)$ its largest singular values. Finally, if $\x$ is a Gaussian random vector with mean $\m$ and covariance $\Q$, we write $\x \sim \mathcal{N}\left(\m,\,\Q\right)$. Similarly, if $x$ is a uniformly distributed random variable between $x_{\min}$ and $x_{\max}$, we denote it as $x \sim \mathcal{U}_{\left[x_{\min},\,x_{\max}\right]}$.
\section{Problem formulation and background}\label{sec:background}

Unconstrained online optimization consists in solving the following minimization problem
\begin{equation}\label{eq:general-online-optimization}
	\x^*_k = \argmin_{\x \in \RNumb^n} f_k(\x) \qquad k \in \N,
\end{equation}
where $f_k\left(\x\right)$ is a time-varying objective function. 

The objective of this paper is to develop control-inspired algorithms to solve a class of online optimization problems in which $f_k\left(\x\right)$ has the following very specific form
\begin{equation}\label{eq:quadratic-cost}
f_k\left(\x\right) =  \frac{1}{2}\left(\x-\c_k\right)^\top\A\left(\x-\c_k\right)+d_k,
\end{equation}
where $\A\in\RNumb^{n\times n}$ is a fixed positive definite matrix and $\c_k\in\RNumb^n$ and $d_k\in\RNumb$ are time varying. It is clear that the optimization problem is independent of $d_k$ (that can hence be fixed to zero) and admits a straightforward solution that is $\x_k^*=\c_k$. The proposed problem makes sense because the only information that is assumed to be available of  $f_k\left(\x\right) $ is the one coming from the evaluation of its gradient in a sequence of points $\x_k$ from which we need to estimate $\c_k$. This is exactly what happens in the standard gradient descent methods that we aim to extend. As in the gradient descent algorithms we assume that the matrix $\A$ is unknown, but that its minimum and maximum eigenvalues are known.

\begin{assum}\label{as:hessianA}
We assume that there exists positive constants $\lambda_{\min},\lambda_{\max}$ such that
\begin{equation}\label{eq:estA}
    \lambda_{\min}\I_n \preceq \A \preceq \lambda_{\max}\I_n,
\end{equation}
We assume moreover that $\lambda_{\min},\lambda_{\max}$ are known.
\end{assum}

\begin{rem}
The online optimization problem \eqref{eq:general-online-optimization} with $f_k\left(\x\right)$ given in \eqref{eq:quadratic-cost} can be seen as the approximation of an optimization problem with a more general cost function that is
\begin{equation}\label{eq:timeTrans}
f_k\left(\x\right) = f\left(\x -\c_k\right)+d_k.
\end{equation}
where $f: \RNumb^n \to \RNumb$ is a fixed smooth strongly convex function with minimum in $\x=0$ and $\c_k\in\RNumb^n$, $d_k\in\RNumb$.

Indeed, by taking the second-order Taylor expansion of $f(\x)$ around $\x=0$
\begin{equation}\label{eq:fapprox}
f(\x)\simeq \frac{1}{2}\x^* \A \x+f(0)
\end{equation}
in which $\A$ is two times the Hessian of $f(\x)$ evaluated at $\x=0$ and the linear term is missing since $f(\x)$ attains the minimum in $\x=0$, then $f(\x_k)$ is well approximated by the quadratic function in \eqref{eq:quadratic-cost}, if we are able to keep the $\x_k$ closed to $\c_k$. This sounds similar to what we do when we approximate a nonlinear model by its linearization, which is reasonable only in case we are able to keep the state closed to a desired value of the state around which we perform the linearization.

Similarly, we can see the proposed online optimization with cost function as in \eqref{eq:quadratic-cost} as the approximation of a general time varying smooth strongly convex cost function $f_k(\x)$ under the assumption that the Hessian around its minimum $\c_k$ is time-invariant (or varies slowly).
\end{rem}

As mentioned above, in the following, we will propose online optimization algorithms providing an estimate $\x_k$ of the minimum $\c_k$ using only the information coming from the gradient of the cost function in $\x_k$. 
Beside being based on the knowledge of of the constants $\lambda_{\min}, \lambda_{\max}$, the design of these algorithms will be based also on a model generating the signal $\c_k$. 
In particular in this paper we make the following assumption regarding the evolution of $\c_k$.

\begin{assum}[Model of $\c_k$]\label{as:modelck}
We assume that the entries $[\c_k]_i$ of  $\c_k$ are generated by the following dynamical system
\begin{equation}\label{eq:linearModels}
\Sigma\,:\begin{cases}
\xiB^{(i)}_{k+1} &= \F\xiB^{(i)}_k + \G [\w_k]_i \\
[\c_k]_i &= \H \xiB^{(i)}_k + j [\w_k]_i
\end{cases}
\end{equation}
where $\xiB^{(i)}_k\in\RNumb^{m}$ is the state vector, $\F\in\RNumb^{m\times m}$, $\G\in\RNumb^{m\times 1}$, $\H\in\RNumb^{1\times m}$, $j\in\RNumb$ and where $\w_k$ is Gaussian white noise with $\w_k\sim\NDist\left(\zero,\,\sigma^2 \I_n\right)$  for all $k$. In other words, $[\w_k]_i$, that are entries of $\w_k$, are independent, zero mean and variance $\sigma^2$ Gaussian white noises.
\end{assum}

We can also say that $\c_k$ evolves according to the extended linear dynamical system
\begin{equation}\label{eq:linearModelsEst}
\Sigma_{ext}\,:\begin{cases}
\xiB_{k+1}^{ext} &= \Fext\xiB_k^{ext} + \Gext \w_k \\
\c_k &= \Hext \xiB_k^{ext} + \Jext \w_k
\end{cases}
\end{equation}
where 
$$\xiB_{k}^{ext}=\left[\begin{array}{c}\xiB^{(1)}_k\\ \vdots\\ \xiB^{(n)}_k\end{array}\right]\in\RNumb^{mn}$$
and $\Fext = \I_n \otimes \F$, $\Gext =\I_n \otimes  \G $, $\Hext =\I_n \otimes \H  $, $\Jext = j \I_n$.

We emphasize that, while we are restricting to the case in which the entries of $\c_k$ are all described by the same model, on the other hand this assumption does not imply that they have exactly the same evolution. Indeed, they are generated by $n$ independent identical SISO systems that act component-wise. We assume that the parameters $\F,\G,\H$ and $j$ are known, while the initial state $\xiB_0^{ext}$ is unknown.

\begin{rem}
The model of $\c_k$ introduced in the previous assumption is equivalent to the existence of a scalar rational transfer function $h(z)$ such that $\c_k$ is the output of a system with transfer matrix $h(z)\I_n$ driven by the Gaussian white noise $\w_k$, where we have that $h\left(z\right) = \H\left(z\I-\F\right)^{-1}\G + j$.
\end{rem}

The model proposed in Assumption  \ref{as:modelck} has two particular cases that are worth to be mentioned:
\begin{enumerate}
\item In case we assume that $\sigma^2=0$ and that $\F$ is marginally unstable, then the entries of $\c_k$ evolve according to the free response of the system $\Sigma$. This case has been treated in \cite{bastianello_internal_2022} and solved using robust control techniques combined with the internal model principle.
\item In case we assume that $\sigma^2>0$ and that $\F$ is stable, then the entries of $\c_k$ are mutually independent Gaussian stationary random processes with the same rational spectra. 
\end{enumerate}

In the following sections, we introduce control-inspired algorithms designed to estimate the solution $\c_k$ of~\eqref{eq:general-online-optimization}.

\section{Algorithms}\label{sec:algorithm}

In this section, we demonstrate how control techniques can be used for designing online optimization algorithms. Specifically, various control concepts can be effectively adapted to address the different challenges posed by online optimization problems, depending on the specific characteristics of the problem at hand.

We propose two algorithms: one inspired by Kalman filtering and the other based on robust $\mathcal{H}_\infty$ control. While the first one can be applied to general models of the signal $\c_k$, the latter is in principle applicable only when the system~\eqref{eq:linearModels} is stable.
For systems that are unstable, the robust $\mathcal{H}_\infty$ controller must be combined with an additional controller to address the instability. 

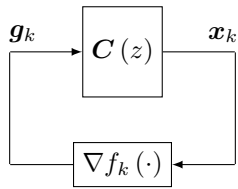
\begin{figure}[!ht]
\centering
\begin{tikzpicture}[auto, >=latex]
	\node[draw,minimum height = 1.2cm] (C) at (0,0.75) {$\C\left(z\right)
		$};
	
	\node[draw=black] (grad) at (0,-0.75) {$\gradf\left(\cdot\right)$};
	
    \node[dot] (d1a) at ($(C)-(1.5,0)$) {};

    \node[dot] (d3a) at ($(C)+(1.5,0)$) {};

    \node[dot] (d2d) at ($(grad)-(-1.5,0)$) {};

    \node[dot] (d4d) at ($(grad)+(1.5,0)$) {};

    \node[dot] (d2a) at (d1a |- d2d) {};
    \node[dot] (d4a) at (d3a |- d4d) {};
	\node[dot] (d2) at (0,0) {};
	\node[dot] (d3) at (0,0) {};
	\node[dot] (d4) at (0,0) {};
	\draw (d1a) edge[->] node[above left] {$\g_k$} (C);

    \draw (d3a) edge[-] node[above right] {$\mathbold{x}_k$} (C);

    \path (d3a) edge (d4a)  (d4a) edge[->] (grad);
    
        \path (d1a) edge (d2a)  (d2a) edge (grad);
		
\end{tikzpicture}
\caption{Proposed control scheme for the solution of~\eqref{eq:general-online-optimization}.}
\label{fig:block-diagram}
\end{figure}

We now introduce the problem more formally. Our objective is to design a controller with transfer matrix $\C\left(z\right)$ that takes as input the gradient signal $\g_k \coloneqq \gradf\left(\x_k\right)$, and produces $\x_k$, which represents estimate of the minimizer $\c_k$ in~\eqref{eq:general-online-optimization}. This corresponds to the block diagram presented in Fig. \ref{fig:block-diagram}.

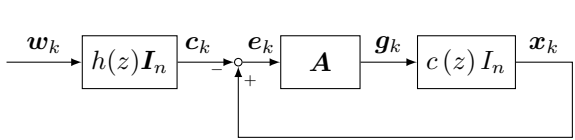
\begin{figure}[!ht]
\centering
\begin{tikzpicture}[auto, >=latex]

\node[input] (input) {};
\node [block, right = of input] (Hwc) {$h(z)\I_n$};
\node [sum, right = 0.75 cm of Hwc] (sum1) {};
\node[input,above = 0.25 cm of sum1] (inputxStar) {};
\node [block, right = 0.5 of sum1] (A) {$\A$};
\node [block, right = 0.75 of A] (controller) {$c\left(z\right)I_n$};
\coordinate[right = 0.75 cm of controller ] (tmp);
\coordinate[below = 1 cm of tmp ] (tmp1);
\draw [draw,->] (controller) -- node[above] {$\x_k$}(tmp) -- (tmp1) -| node [xshift = 3.75mm,yshift = 7.5mm] {\tiny $+$} (sum1);
\draw [draw,->] (sum1) -- node[above] {$\e_k$}(A);
\draw [draw,->] (A) -- node[above] {$\g_{k}$}  (controller);
\draw [draw,->] (input) -- node[above] {$\w_{k}$}  (Hwc);
\draw [draw,->] (Hwc) -- node[above,xshift = -0.1cm] {$\c_k$} node [yshift = -3mm, xshift = 1.5mm] {\tiny $-$} (sum1);
\end{tikzpicture}
\caption{Control scheme of Fig.~\ref{fig:block-diagram} when $ \C(z) = c(z) \I_n $ and condition~\eqref{eq:gradE} holds.
}
\label{fig:block-diagram-caseA}
\end{figure}

In general, the controller $\C\left(z\right)$ is a general rational strictly proper transfer matrix.
In this work, since in our set-up the matrix $\A$ is not known (see Assumption~\ref{as:hessianA}), it is convenient to impose that $\C\left(z\right) = c\left(z\right)\I_n$, where $c\left(z\right)$  is a scalar strictly proper rational transfer function.
Under this assumption, and since in our case
\begin{equation}\label{eq:gradE}
    \g_k \coloneqq  \gradf\left(\x_k\right) = \A\underbrace{\left(\x_k-\c_k\right)}_{ \e_k\coloneqq}
\end{equation}
we can transform the feedback scheme in Fig.~\ref{fig:block-diagram} into the one shown in Fig.~\ref{fig:block-diagram-caseA}

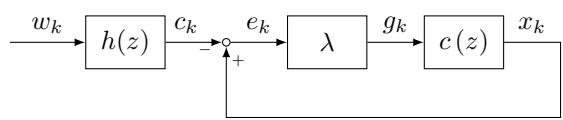
\begin{figure}[!ht]
\centering
\begin{tikzpicture}[auto, >=latex]
\node[input] (input) {};
\node [block, right = of input] (Hwc) {$h(z)$};
\node [sum, right = 0.75 cm of Hwc] (sum1) {};
\node[input,above = 0.25 cm of sum1] (inputxStar) {};
\node [block, right = 0.75 of sum1] (A) {$\lambda$};
\node [block, right = 0.75 of A] (controller) {$c\left(z\right)$};
\coordinate[right = 0.75 cm of controller ] (tmp);
\coordinate[below = 1 cm of tmp ] (tmp1);
\draw [draw,->] (controller) -- node[above] {$x_k$}(tmp) -- (tmp1) -| node [xshift = 3.75mm,yshift = 7.5mm] {\tiny $+$} (sum1);
\draw [draw,->] (sum1) -- node[above] {$e_k$}(A);
\draw [draw,->] (A) -- node[above] {$g_k$}  (controller);
\draw [draw,->] (input) -- node[above] {$w_k$}  (Hwc);
\draw [draw,->] (Hwc) -- node[above,xshift = -0.1cm] {$c_k$} node [yshift = -3mm, xshift = 1.5mm] {\tiny $-$} (sum1);
\end{tikzpicture}
\caption{Block scheme representation of the transfer function~\eqref{eq:scalarTF}.}
\label{fig:block-diagram-caseScalar}
\end{figure}

Due to Assumption~\ref{as:modelck}, the minimizer $\c_k$ is driven by white noise. Therefore, the error $\e_k := \x_k - \c_k$ cannot converge to zero. Instead, we can try to analyze the statistical properties of $\e_k$. Precisely, we consider the following cost function
\begin{equation}\label{eq:costJ1}
    J \coloneqq \lim_{k\to +\infty}\mathbb{E}\left[\e_k^{\top}\e_k\right],
\end{equation}
and our objective will be to find the controller $ c(z) $ that minimizes this cost.
This is also the $ \mathcal{H}_2 $- norm of the transfer matrix $\H_{\w\e}\left(z\right)$ that is the transfer matrix from the noise input $\w$ to the error output $\e$ in the control scheme shown in Fig.~\ref{fig:block-diagram-caseA}, namely we can write
\begin{equation}\label{eq:costJ2}
    J = \norm{\H_{\w\e}\left(z\right)}^2_{2}.
\end{equation}
Observe that
\begin{equation}\label{eq:errorTM}
    \H_{\w\e}\left(z\right) = -\left(\I_n - c\left(z\right)\A \right)^{-1} h\left(z\right).
\end{equation}
Ideally, we would like to combine~\eqref{eq:costJ2} and~\eqref{eq:errorTM} to determine an optimal controller $c\left(z\right)$. However, we do not know $\A$, and thus must account for this uncertainty. The following calculations are useful for this purpose, as they demonstrate that the minimization of~\eqref{eq:costJ2} depends only on the eigenvalues of $\A$ and it is independent of its eigenvectors. Indeed, by~\eqref{eq:costJ2}, we have
\begin{equation}\label{eq:intComp}
    J = \frac{1}{2\pi}\int_{-\pi}^{\pi}\tr\left[\H_{\w\e}^{*}\left(e^{i\theta}\right)\H_{\w\e}\left(e^{i\theta}\right)\right] \,d\theta.
\end{equation}
Since $\A$ is symmetric, we can always consider its diagonalization $ \A = \V \LB\V^{\top} $, where $\LB$ is real and diagonal, and $\V$ is a real orthogonal matrix. This decomposition allows us to obtain that
\begin{equation}\label{eq:Hew}
    \H_{\w\e}\left(z\right) = -\V\left(\I_n - c\left(z\right)\LB\right)^{-1} \V^{*}h\left(z\right),
\end{equation}
which yields 
\begin{equation}\label{eq:traceIntComp}
    \tr\left[\H_{\w\e}^{*}\left(z\right)\H_{\w\e}\left(z\right)\right] =\sum_{i=1}^{n}\frac{ \abs{h\left(z\right)}^2}{\abs{1-\lambda_ic\left(z\right)}^2}.
\end{equation}
Substituting~\eqref{eq:traceIntComp} into Equation~\eqref{eq:intComp} yields
\begin{equation}\label{eq:Jcost}
\begin{aligned}
    J &= \sum_{i=1}^{n}\overbrace{\frac{1}{2\pi}\int_{-\pi}^{\pi}\frac{\abs{h\left(e^{i\theta}\right)}^2}{\abs{1-\lambda_ic\left(e^{i\theta}\right)}^2}\,d\theta}^{J\left(\lambda_i,\,c(z)\right)\coloneqq}\\
    &=\sum_{i=1}^{n}J\left(\lambda_i,\,c(z)\right).
\end{aligned}
\end{equation}
The problem we aim to solve is to determine the controller $c(z)$ able to minimize the cost $J$.  
Each term $J\left(\lambda_i,\,c(z)\right)$ in the sum providing the cost $J$ is itself the $ \mathcal{H}_2$- norm of the transfer function 
\begin{equation}\label{eq:scalarTF}
    w_\lambda(z):=-\frac{h(z)}{1-\lambda c(z)},
    \end{equation}
    with $\lambda = \lambda_i$.
Notice that $w_\lambda(z)$ is the transfer function from the input $w_k$ to the output $e_k$ in the block diagram in Fig. \ref{fig:block-diagram-caseScalar}.


It is important to note that the exact values of the eigenvalues $\lambda_i$ are not known. However, based on condition~\eqref{eq:estA}, we do know that they satisfy $\lambda_{\min}\leq\lambda_i\leq\lambda_{\max}$. To address this uncertainty, we propose two robust approaches: the first is a Kalman-inspired method, and the second is a robust $\mathcal{H}_\infty$-inspired method.


\subsection{Kalman-inspired algorithm}

In order to explain the ratio of this method we get inspiration to what is the optimal choice in the limit case in which $n=1$ so that the cost $J$ is composed only of one term $J=J(\lambda,c(z))$ and hence the optimal controller has to be designed so to minimize such a cost. 
We see that the optimal controller is $c_\lambda(z)$ that is realized by the following state space iteration
\begin{equation}\label{eq:estimator1}
\begin{aligned}
    \xiH_{k+1} &= \F \xiH_k - \lambda^{-1}\K g_k \\
    x_k &= \H \xiH_k,
\end{aligned}
\end{equation}
where $\xiH_k\in\RNumb^m$, $\F$ and $\H$ are the same matrices as in~\eqref{eq:linearModels}, and $\K$ is the steady-state Kalman filter gain~\cite[Theorem 8.5.3]{van_control_2021} that is given by the formula
\begin{equation}\label{eq:KFgain}
    \K = \left(\F\P\H^\top + \G\sigma^2j\right)\left(\H\P\H^\top+j^2\right)^{-1},
\end{equation}
where $\P$ is the solution of the following algebraic Riccati equation
\begin{multline*}
    \P = \F\P\F^\top+\sigma^2\G\G^\top-\left(\F\P\H^\top + \G\sigma^2j\right)\\\cdot\left(\H\P\H^\top+j^2\right)^{-1}\left(\F\P\H^\top + \G\sigma^2j\right)^\top.
\end{multline*}
Indeed, notice that the previous system equations combined with the signal relations displayed in the block diagram in Fig.~\ref{fig:block-diagram-caseScalar}, yields 
\begin{equation}\label{eq:kfeq}
\begin{aligned}
    \xiH_{k+1} &= \F \xiH_k +\K(c_k- \H \xiH_k)\\
    x_k &= \H \xiH_k,
\end{aligned}
\end{equation}
that are the Kalman filter equations associated with the signal generated by the stochastic model
\begin{equation}\label{eq:stocmodel}
\begin{cases}
\xiB_{k+1} &= \F\xiB_k + \G w_k \\
c_k &= \H \xiB_k + j w_k
\end{cases}
\end{equation}
where $\F,\G,\H,j$ are the same as in Assumption \ref{as:modelck} and $w_k$ is a zero mean and variance $\sigma^2$ Gaussian white noise.
Indeed, $e_k$ is the innovation of this Kalman filter whose asymptotic covariance is minimum possible by the Kalman filtering optimality~\cite[Proposition 8.6.1]{van_control_2021}. The controller $c_\lambda(z)$ continues to remain optimal if $n>1$ and the eigenvalues $\lambda_i$ are all equal to $\lambda$. 
If we assume that the eigenvalues $\lambda_i$ are not equal but are only closed to $\lambda$, then we can expect that a perturbation of the controller $c_\lambda(z)$ will be closed to be optimal.

For this reason we propose here a sub-optimal controller keeping the structure of $c_\lambda(z)$ but in which the parameter $\lambda$ is tuned so to minimize the cost $J$, namely we want to solve the following optimization problem 
 \begin{equation}\label{eq:optH2}
     \lambda^* = \argmin_{\lambda > 0 } \sum_{i=1}^{n}J\left(\lambda_i,\,c_\lambda(z)\right).
 \end{equation}
 
 Notice that, if we apply a controller $c_{\mu}(z)$ with $\mu\not=\lambda$ in the interconnection in Fig. \ref{fig:block-diagram-caseScalar}, the resulting dynamics does not yield to the Kalman filter but instead to the iteration
\begin{equation}\label{eq:kfeqapprox}
\begin{aligned}
    \xiH_{k+1} &= \F \xiH_k +\frac{\lambda}{\mu}\K(c_k- \H \xiH_k)\\
    x_k &= \H \xiH_k,
\end{aligned}
\end{equation}
The resulting cost $J(\lambda,c_\mu(z))$ will be a function of $\lambda/\mu$, namely 
$$J(\lambda,c_\mu(z))=J(\lambda/\mu)$$
We know that the minimum of $J(a)$ is attained in $a=1$. We approximate this function by a second-order expansion around $a = 1$
\begin{equation}
    J\left(a\right) \approx J\left(1\right) + \underbrace{\left.\frac{\partial J\left(a\right)}{\partial a}\right\rvert_{a = 1}}_{ = 0}\left(a-1\right) + \underbrace{\left.\frac{\partial^2 J\left(a\right)}{\partial a^2}\right\rvert_{a = 1}}_{\beta}\left(a-1\right)^2,
\end{equation}
where the linear term is zero due to the optimality of $a = 1$. This approximation holds if the fractions $a=\lambda/\mu$ is close to $1$. 

If we apply the controller $c_{\mu}(z)$ in the scheme in Fig. \ref{fig:block-diagram-caseA}, under the assumption that all the $\lambda_i$'s are closed to $\mu$, we can rewrite~\eqref{eq:Jcost} as
\begin{equation}\label{eq:costJApprox}
\begin{aligned}
    J &\approx \sum_{i=1}^n \beta \left(\frac{\lambda_i}{\mu}-1\right)^2 + \text{const.}\\
    &= \beta\left(\frac{\sum_{i=1}^n\lambda_i^2}{\mu^2}-2\frac{\sum_{i=1}^n\lambda_i}{\mu}\right)+ \text{const.}
\end{aligned}
\end{equation}
which is quadratic in $1/\mu$. By this approximation the solution of the optimization problem~\eqref{eq:optH2} is
\begin{equation}\label{eq:muOpt}
    \mu^* = \frac{\sum_{i=1}^n\lambda_i^2}{\sum_{i=1}^n\lambda_i}.
\end{equation}
This optimal value $\mu^*$ can be computed only when the eigenvalues of $\A$ are known, and this is not typically the case. If we know only a probability distribution $\phi\left(\lambda\right)$ of the eigenvalues, then we can take  
\begin{equation}\label{eq:muOptExp}
    \mu^*= \frac{\int\lambda^2\phi(\lambda)d\lambda}{\int\lambda\phi(\lambda)d\lambda}\ .
\end{equation}
 In practice, equation~\eqref{eq:muOptExp} can be applied by incorporating prior knowledge of the distribution of the eigenvalues of $\A$ into $\phi\left(\lambda\right)$. 

Without additional information, we can assume that $\lambda_i$ are uniformly distributed, namely $\lambda\sim\mathcal{U}_{\left[\lambda_{\min},\,\lambda_{\max}\right]}$. Under this assumption, from~\eqref{eq:muOptExp} we obtain
\begin{align}\label{eq:muOptFinalUniform}
    \mu^* &= \frac{\int_{\lambda_{\min}}^{\lambda_{\max}}\lambda^2\phi\left(\lambda\right)\,d\lambda}{\int_{\lambda_{\min}}^{\lambda_{\max}}\lambda\phi\left(\lambda\right)\,d\lambda}= \frac{2}{3}\frac{\lambda_{\max}^2 + \lambda_{\max}\lambda_{\min} + \lambda_{\min}^2}{\lambda_{\max} + \lambda_{\min}}.
\end{align}
Using the optimal $\mu^*$, we can implement the controller $\C\left(z\right) = c_{\mu^*}\left(z\right)\I_n$ in Fig.~\ref{fig:block-diagram} that is
\begin{equation}\label{eq:controllerExtEst}
\begin{aligned}
    \xiH_{k+1}^{ext} &= \Fext\xiH_k^{ext} -\frac{1}{\mu^*}\Kext \gradf\left(\x_k\right)\\
    \x_k &= \Hext\xiH_k^{ext},
\end{aligned}
\end{equation}
where $\Fext$, and $\Hext$ are the same matrices of~\eqref{eq:linearModelsEst}, $\Kext = \I_n \otimes \K$, with $\K$ being the Kalman filter gain in \eqref{eq:KFgain}, and $\xiH_k^{ext}\in\RNumb^{nm}$ is the controller state.
\begin{rem}
    It is important to note that the representation in Fig.~\ref{fig:block-diagram-caseScalar} is primarily intended for analysis purposes. This representation allows us to interpret the scheme in Fig.~\ref{fig:block-diagram-caseA} as a collection of $n$ independent SISO systems, each corresponding to Fig.~\ref{fig:block-diagram-caseScalar}. For practical implementation, we directly consider the original scheme in Fig.~\ref{fig:block-diagram-caseA}, where the controller $\C\left(z\right)=c\left(z\right)\I_n$ takes $\g_k$ as input and returns $\x_k$ as output, as shown in Fig.~\ref{fig:block-diagram}. The signals $g_k$ and $x_k$ in Fig.~\ref{fig:block-diagram-caseScalar} are auxiliary signals used solely for analysis and design purposes, and knowledge of the decomposition $\A = \V\LB\V^{\top}$ is not required for the implementation.
\end{rem}
\subsection{Robust $\mathcal{H}_\infty$-inspired algorithm}
    While our primary objective is to minimize the cost in~\eqref{eq:costJ2}, due to its relationship with the error statistics in~\eqref{eq:costJ1}, it is important to note that the $\mathcal{H}_\infty$-norm serves as an upper bound to the $\mathcal{H}_2$-norm~\cite{debruyne_linear_1995}. This bound provides a valuable tool for deriving a suboptimal solution to the minimization of $J$ in~\eqref{eq:costJ2} using robust $\mathcal{H}_{\infty}$-control tools. Specifically, we know that
\begin{equation}
    J = \norm{\H_{\w\e}\left(z\right)}_{2} \leq \norm{\H_{\w\e}\left(z\right)}_{\infty}.
\end{equation}
Using the definition of the $\mathcal{H}_{\infty}$-norm, we have
\begin{equation}\label{eq:HinfNorm}
    \norm{\H_{\w\e}\left(z\right)}_{\infty} \coloneqq \sup_{\theta \in \left[-\pi,\, \pi\right)}\sigma_{\max}\left(\H_{\w\e}\left(e^{i\theta}\right)\right),
\end{equation}
where $\H_{\w\e}\left(z\right)$ is defined as in~\eqref{eq:Hew}. Since the singular values are invariant under unitary transformations, we can rewrite this as
\begin{align*}
    &\norm{\H_{\w\e}\left(z\right)}_{\infty} \\
    &\qquad =\sup_{\theta \in \left[-\pi,\, \pi\right)}\sigma_{\max}\left(-\left(\I_n - c\left(e^{i\theta}\right)\LB\right)^{-1}h\left(e^{i\theta}\right)\right)\\
    &\qquad= \max_{i=1,\ldots,n}\sup_{\theta \in \left[-\pi,\, \pi\right)}\frac{\abs{h\left(e^{i\theta}\right)}}{\abs{1 - \lambda_{i}c\left(e^{i\theta}\right)}}.
\end{align*}

Following a similar approach to the previous formulations, we now define the robust control problem, that is finding a controller $c\left(z\right)$ that minimizes the cost function
\begin{align}
    \Jhat &= \sup_{\lambda\in[\lambda_{\min},\lambda_{\max}]}\sup_{\theta \in \left[-\pi,\, \pi\right)}\frac{\abs{h\left(e^{i\theta}\right)}}{\abs{1 - \lambda c\left(e^{i\theta}\right)}}\\
    &=\sup_{\lambda\in[\lambda_{\min},\lambda_{\max}]}\norm{w_{\lambda}\left(z\right)}_\infty.
\end{align}
where $w_{\lambda}\left(z\right)$ is defined in~\eqref{eq:scalarTF}\footnote{In practice, this optimization problem is solved by iteratively decreasing a threshold $\gamma > 0$ and finding a controller $c\left(z\right)$ such that $\norm{w_{\lambda}\left(z\right)}_\infty \leq \gamma$ for $\lambda_{\min}\leq\lambda\leq \lambda_{\max}$, until no such controller can be found.}.
In this way we treat the eigenvalues as an uncertainty.
This formulation of our problem is a standard $\mathcal{H}_\infty$-robust optimal control problem, which can be solved using the $\mu$-synthesis methodology, a widely used  strategy for addressing robust $\mathcal{H}_{\infty}$-optimal control problems, as described in~\cite[Section 10]{zhou_robust_1998}. 

In this subsection, while we do not present the complete algorithm for performing $\mu$-synthesis as detailed in~\cite[Section 10.4]{zhou_robust_1998}, we describe how our setup can be adapted for this application. Notice that, when $\F$ in~\eqref{eq:linearModels} is stable, $\mu$-synthesis is directly applicable. However, if $\F$ is unstable (for instance only marginally stable), a two-stage design procedure is needed. This is due to the fact that the unstable poles of $h\left(z\right)$ act as uncontrollable poles in the feedback system shown in Fig.~\ref{fig:block-diagram-caseScalar}, violating the assumptions of $\mathcal{H}_\infty$ theory~\cite[Sections 14.7-8]{zhou_robust_1998}.

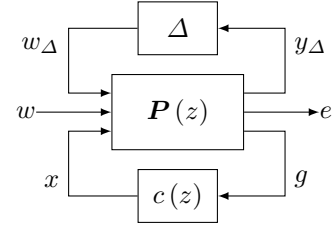
\begin{figure}[!ht]
\centering
\begin{tikzpicture}[auto, >=latex]
\node[input] (w){$w$};
\node[input, above = 0.15cm of w] (wDelta){};
\node[input, below = 0.15cm of w] (u){};
\node [block, minimum width = 1.75cm,minimum height = 1 cm, right = of w] (P) {$\P\left(z\right)$};
\node [block, above = 0.25cm of P] (Delta) {$\Delta$};
\node [block, below = 0.25cm  of P] (ctrl) {$c\left(z\right)$};
\node[output,right = of P] (e){$e$};
\node[output, above = 0.15cm of e] (eDelta){};
\node[output, below = 0.15cm of e] (y){};
\draw [draw,->] (w) -- (P);
\draw [draw,->] (P) -- (e);
\coordinate (tmpInputDelta) at ($(P.west)!0.5!(P.north west)$);
\coordinate (tmpInputU) at ($(P.west)!0.5!(P.south west)$);
\coordinate (tmpMiddleInput) at ($(w)!0.5!(P.west)$);
\coordinate (tmpMiddleOutput) at ($(e)!0.5!(P.east)$);
\coordinate (tmpOutputDelta) at ($(P.east)!0.5!(P.north east)$);
\coordinate (tmpOutputU) at ($(P.east)!0.5!(P.south east)$);

\coordinate (tmpMiddleInputUp) at (tmpInputDelta-|tmpMiddleInput);
\coordinate (tmpMiddleOutputUp) at (tmpOutputDelta-|tmpMiddleOutput);

\coordinate (tmpMiddleInputDown) at (tmpInputU-|tmpMiddleInput);
\coordinate (tmpMiddleOutputDown) at (tmpOutputU-|tmpMiddleOutput);

\draw [draw,->] (tmpOutputDelta) -- (tmpMiddleOutputUp)  |- node[below right]{$y_{\Delta}$} (Delta);
\draw [draw,->] (tmpOutputU) -- (tmpMiddleOutputDown) |- node[above right]{$g$} (ctrl);
\draw [draw,->] (Delta) -| node[below left]{$w_{\Delta}$} (tmpMiddleInputUp) -- (tmpInputDelta) ;
\draw [draw,->] (ctrl)-|node[above left]{$x$}  (tmpMiddleInputDown) -- (tmpInputU);
\end{tikzpicture}
\caption{General framework for $\mu$-synthesis.}
\label{fig:block-diagram-generalFramework}
\end{figure}
We are now in a position to align our framework with the general approach described in~\cite{zhou_robust_1998}. This alignment is achieved by reshaping the control scheme in Fig.~\ref{fig:block-diagram-caseScalar} to match the standard scheme shown in Fig.~\ref{fig:block-diagram-generalFramework}.
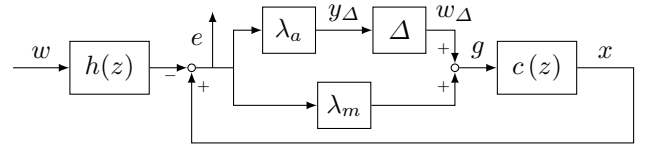
\begin{figure}[!ht]
\centering
\begin{tikzpicture}[auto, >=latex]
\node[input] (input) {};
\node [block, right =0.75  of input] (Hwc) {$h(z)$};
\node [sum, right = 0.5 of Hwc] (sum1) {};
\node[input,above = 0.25 cm of sum1] (inputxStar) {};
\node [block, right = 4cm of sum1] (controller) {$c\left(z\right)$};
\node [sum, left = 0.5 cm of controller] (sum2) {};

\coordinate[right = 0.5 cm of sum1 ] (tmpLIn);
\coordinate (tmpOutput) at ($(sum1)!0.5!(tmpLIn)$);
\node[input,above = 0.75 cm of tmpOutput] (e) {};
\draw [draw,->] (tmpOutput) -- node  {$e$}(e);

\coordinate (tmpLInMiddleOut) at ($(tmpLIn)!0.5!(sum2)$);
\node [blockSmall, below = 0.5cm of tmpLInMiddleOut,anchor = center] (lm) {$\lambda_m$};
\coordinate[above = 0.5cm of tmpLInMiddleOut ] (tmpMiddleUp);
\node [blockSmall, right = 0.375cm of tmpMiddleUp] (delta) {$\Delta$};
\node [blockSmall, left = 0.375cm of tmpMiddleUp] (la) {$\lambda_a$};
\coordinate[right = 0.75 cm of controller ] (tmp);
\coordinate[below = 1 cm of tmp ] (tmp1);
\draw [draw,->] (controller) -- node[above] {$x$}(tmp) -- (tmp1) -| node [xshift = 3.75mm,yshift = 7.5mm] {\tiny $+$} (sum1);
\draw [draw] (sum1) -- (tmpLIn);
\draw [draw,->] (tmpLIn) |- (la);
\draw [draw,->] (la) -- node [above, yshift = 0mm, xshift = 0mm] {$y_\Delta$} (delta);
\draw [draw,->] (tmpLIn) |- (lm);
\draw [draw,->] (lm) -| (sum2) node [yshift = -2.5mm, xshift = -1.5mm] {\tiny $+$};
\draw [draw,->] (delta) -| node [above, yshift = 0mm, xshift = 0mm] {$w_\Delta$} (sum2) node [yshift = 2.5mm, xshift = -1.5mm] {\tiny $+$};
\draw [draw,->] (sum2) -- node [above, yshift = 0mm, xshift = 0mm] {$g$} (controller);
\draw [draw,->] (input) -- node[above] {$w$}  (Hwc);
\draw [draw,->] (Hwc) -- node[above,xshift = -0.1cm] {} node [yshift = -3mm, xshift = 0.25mm] {\tiny $-$} (sum1);
\end{tikzpicture}
\caption{General framework shown in Fig.~\ref{fig:block-diagram-generalFramework}, tailored to the specific control scheme depicted in Fig.~\ref{fig:block-diagram-caseScalar}. }
\label{fig:block-diagram-generalFramework2}
\end{figure}
Specifically, we rewrite the scheme in Fig.~\ref{fig:block-diagram-caseScalar} as depicted in Fig.~\ref{fig:block-diagram-generalFramework2}, where $\Delta \in \left[-1, 1\right]$, and $w_{\Delta}$ and $y_{\Delta}$ are auxiliary signals used to capture system uncertainties~\cite[Section 3.4]{scherer_theory_2001}. Here, we define $\lambda_m = \left(\lambda_{\max} - \lambda_{\min}\right)/2$ and $\lambda_a = \left(\lambda_{\max} + \lambda_{\min}\right)/2$, so that $\lambda = \lambda_m + \lambda_a \Delta$. The entries of the transfer matrix $\P\left(z\right)$ are then derived by a direct computation from the scheme in Fig.~\ref{fig:block-diagram-generalFramework2}.

\subsubsection{$\F$ stable}
In this case, the assumptions of robust $\mathcal{H}_{\infty}$-control design are satisfied, enabling the direct application of $\mu$-synthesis. This process yields a scalar controller $c\left(z\right)$, which we require to be strictly proper to avoid algebraic loops in the control scheme. Given a state-space realization $\left(\Fhat,\,\Ghat,\,\Hhat\right)$ of $c\left(z\right)$, we define the controller $\C\left(z\right) = c\left(z\right)\I_n$ as shown in Fig.~\ref{fig:block-diagram} by the equations
\begin{equation}\label{eq:controllerExtInt}
\begin{aligned}
    \xiH_{k+1}^{ext} &= \Fhatext\xiH_k^{ext}+ \Ghatext \gradf\left(\x_k\right)\\
    \x_k &= \Hhatext\xiH_k^{ext},
\end{aligned}
\end{equation}
where $\Fhatext = \I_n \otimes \Fhat$, $\Ghatext = \I_n \otimes \Ghat$, and $ \Hhatext = \I_n \otimes \Hhat$. 

\subsubsection{$\F$ unstable}
Since the poles of $h\left(z\right)$ are uncontrollable in the scheme shown in Fig.~\ref{fig:block-diagram-generalFramework}, if $\F$ is unstable, then the direct application of $\mathcal{H}_\infty$-robust control design is not possible. However, by reformulating the problem and introducing a precompensator in the control loop, it is possible to derive a controller $c\left(z\right)$ that satisfies the desired properties.

Suppose we factorize $h\left(z\right)$ as follows:
\begin{equation}\label{eq:hduds}
    h\left(z\right) = \frac{n\left(z\right)}{d_u\left(z\right)d_s\left(z\right)}
\end{equation}
where the numerator $n\left(z\right)$ and the denominator $d_u\left(z\right)d_s\left(z\right)$ is factorized so that $d_u\left(z\right)$ contains the unstable poles of $h(z)$, while $d_s\left(z\right)$ contains the stable poles. We assume that $c\left(z\right)$ has the structure
\begin{equation}\label{eq:controllerParUnstable}
    c\left(z\right) =\frac{p\left(z\right)}{d_u\left(z\right)}\cbarS\left(z\right)
\end{equation}
where $p\left(z\right)$ is an arbitrary stable polynomial of the same degree as $d_u\left(z\right)$.
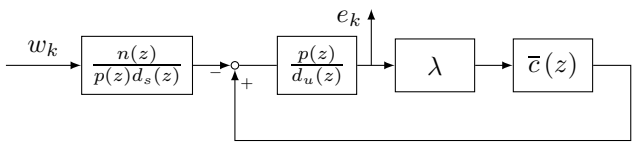
\begin{figure}[!ht]
\centering
\begin{tikzpicture}[auto, >=latex]
\node[input] (input) {};
\node [block, right = of input] (Hwc) {$\frac{n\left(z\right)}{p\left(z\right)d_s\left(z\right)}$};
\node [sum, right = 0.5 cm of Hwc] (sum1) {};
\node[input,above = 0.25 cm of sum1] (inputxStar) {};
\node [block, right = 0.5 of sum1] (precomp) {$\frac{p\left(z\right)}{d_u\left(z\right)}$};
\node [block, right = 0.5 of precomp] (A) {$\lambda$};
\node [block, right = 0.5 of A] (controller) {$\cbarS\left(z\right)$};
\coordinate[right = 0.5 cm of controller ] (tmp);
\coordinate[below = 1 cm of tmp ] (tmp1);
\draw [draw,->] (controller) -- (tmp) -- (tmp1) -| node [xshift = 3.75mm,yshift = 7.5mm] {\tiny $+$} (sum1);

\coordinate (tmpOutput) at ($(precomp.east)!0.375!(A.west)$);
\node[input,above = 0.75 cm of tmpOutput] (e) {};
\draw [draw,->] (tmpOutput) -- node[yshift = 0.25cm]  {$e_k$} (e);
\draw [draw,->] (sum1) -- (precomp) -- (A);
\draw [draw,->] (A) --   (controller);
\draw [draw,->] (input) -- node[above] {$w_k$}  (Hwc);
\draw [draw,->] (Hwc) -- node[below,xshift = -0.1cm] {} node [yshift = -3mm, xshift = 0.5mm] {\tiny $-$} (sum1);
\end{tikzpicture}
\caption{Equivalent control scheme to Fig.~\ref{fig:block-diagram-caseScalar} with $c\left(z\right)$ parametrized as in Equation~\eqref{eq:controllerParUnstable}.}
\label{fig:block-diagram-caseScalarUnstable}
\end{figure} 
It is easily seen that the transfer function from the input $w_k$ to the output $e_k$ in the scheme in Fig.~\ref{fig:block-diagram-caseScalar} is the same as the transfer function from the input $w_k$ to the output $e_k$ in the scheme in Fig. \ref{fig:block-diagram-caseScalarUnstable}. In the new scheme however the uncontrollable poles are all stable. By following the previous procedure it is possible 
to derive a controller $\cbarS\left(z\right)$ via $\mu$-synthesis, which can then be used to obtain $c\left(z\right)$ by~\eqref{eq:controllerParUnstable} and the final controller $\C\left(z\right) = c\left(z\right)\I_n$ whose realization is given in~\eqref{eq:controllerExtInt}.
Notice that this method introduces a degree of freedom in the choice of the polynomial $p(z)$, that, for simplicity, will be set to $p\left(z\right) = z^{n_u} $, where $n_u$ is the degree of $d_u(z)$.



\section{Simulations}\label{sec:simulations}
In this section, we present numerical experiments in order to analyse the performance of the proposed algorithms. Moreover, we compare our algorithms with the standard online version of gradient descent algorithm, that is
\begin{equation}\label{eq:ogd}
\x_{k+1} = \x_{k}-\alpha \gradf\left(\x_k\right)
\end{equation}
with $\alpha$ set to $1/\lambda_{\max}$. Additionally, when $\F$ is unstable, we compare the Kalman-inspired algorithm and the precompensated robust $\mathcal{H}_{\infty}$-inspired algorithm against a modified version of the approach proposed in~\cite{bastianello_internal_2022}, which has been optimized for obtaining the best possible convergence rate.
For all the following simulations, we fix $n = 10$ in~\eqref{eq:quadratic-cost} and $\sigma = 1$ in~\eqref{eq:linearModels}.

\subsection{$\F$ stable}
We first consider a case in which $\F$ is assumed to be stable. We compare the gradient descent algorithm~\eqref{eq:ogd}, the Kalman-inspired algorithm~\eqref{eq:controllerExtEst}, and the robust $\mathcal{H}_\infty$-inspired algorithm~\eqref{eq:controllerExtInt}.

The first experiment illustrates the time evolution of $\norm{\e_k}$ where $\e_k:=\x_k - \c_k$ over time $k$. The simulation uses a randomly generated matrix $\A$. Specifically, we generate $\lambda_i$ uniformly distributed between $\lambda_{\min}$ and $\lambda_{\max}$, we define $\LB$ as a diagonal matrix with diagonal entries $\lambda_i$, and we generate a random orthogonal matrix $\V$, setting $\A=\V\LB\V^*$. The pair $\left(\F,\,\H\right)$ in~\eqref{eq:linearModels} is defined in canonical observable form~\cite{joao_linear_2018} with the characteristic polynomial of $\F$ given by $d_s\left(z\right) = \left(z-p\right)^2$, where $p = 0.975$. We set $\G$ to a vector of all ones. 
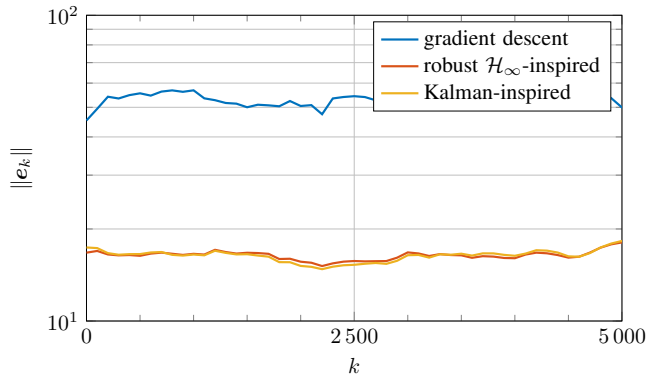
\begin{figure}[!ht]
\centering
\input{./Figures/figTest0.tex}
\caption{Comparison of the time evolution of the moving average of the error norm $\norm{\e_k}$, averaged over $1000$ local points, in a logarithmic scale for the online gradient descent, robust $\mathcal{H}_{\infty}$-inspired~\eqref{eq:controllerExtInt}, and Kalman-inspired~\eqref{eq:controllerExtEst} algorithms over time $k$.}
\label{fig:figTest0}
\end{figure}

The results of this initial numerical experiment are shown in Fig.~\ref{fig:figTest0}. This figure illustrates the performance differences between the proposed algorithms (Kalman-inspired and robust $\mathcal{H}_\infty$-inspired) and the traditional gradient descent method. Notably, the proposed algorithms demonstrate improved performance.

To further evaluate and distinguish the performance of these approaches, we conduct additional experiments and plot the cost $J$ as a function of varying parameters.
\begin{figure}[!ht]
\centering
\input{./Figures/figTest1.tex}
\caption{Square root of the cost $J$ from~\eqref{eq:Jcost} for a stable $\F$, shown on a logarithmic scale as a function of measurement noise level $j$, comparing the performance of online gradient descent, robust $\mathcal{H}_{\infty}$-inspired, and Kalman-inspired algorithms.}%
\label{fig:figTest1}
\end{figure}
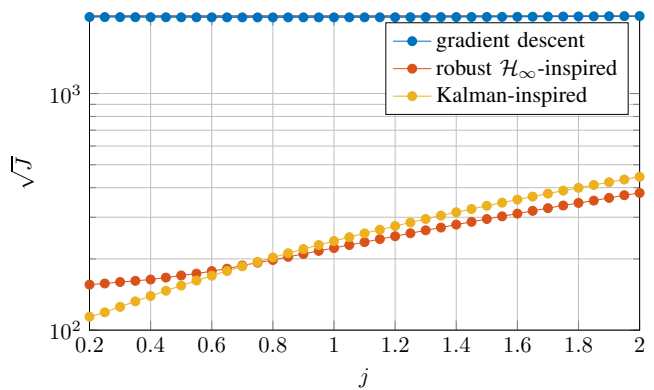
The first of these tests consists keeping the same model generating $\c_k$ as above, except that we vary $j$ in the interval $\left[0.2,\,2\right]$. The results of this test are shown in Fig.~\ref{fig:figTest1}.
\begin{figure}[!ht]
\centering
\input{./Figures/figTest2.tex}
\caption{Square root of the cost $J$ from~\eqref{eq:Jcost} with a stable $\F$, shown on a logarithmic scale as a function of $\lambda_{\max}$, comparing the performance of online gradient descent, robust $\mathcal{H}_{\infty}$-inspired, and Kalman-inspired algorithms.}%
\label{fig:figTest2}
\end{figure}
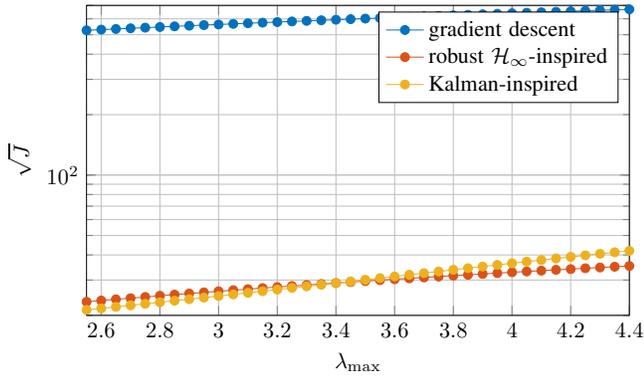%

The second test, depicted in Fig.~\ref{fig:figTest2}, uses the same model, setting $j = 0.2$ and varying $\lambda_{\max}$ within the interval $\left[2.6,\,4.4\right]$.

These two experiments reveal performance differences between the approaches. The comparative analysis of Figs.~\ref{fig:figTest1} and~\ref{fig:figTest2} suggests that while both proposed methods outperform gradient descent, the Kalman-inspired algorithm may be preferable when the system dynamics are well known or less noisy. Conversely, under higher uncertainty, the robust $\mathcal{H}_\infty$-inspired algorithm appears to yield better results.

\subsection{$\F$ unstable}

We now evaluate the algorithms in the scenario where $\F$ has eigenvalues on the unit circle. As before, we conduct tests similar to those presented in Figs.~\ref{fig:figTest1} and~\ref{fig:figTest2}. Though, in these cases, we cannot compare our methods with the gradient descent algorithm because the cost function $J$ diverges, resulting in an infinite value, we compare our approaches against the method proposed in~\cite{bastianello_internal_2022}, which is optimized for convergence rate and follows a design procedure based on the Internal Model Principle (IMP).
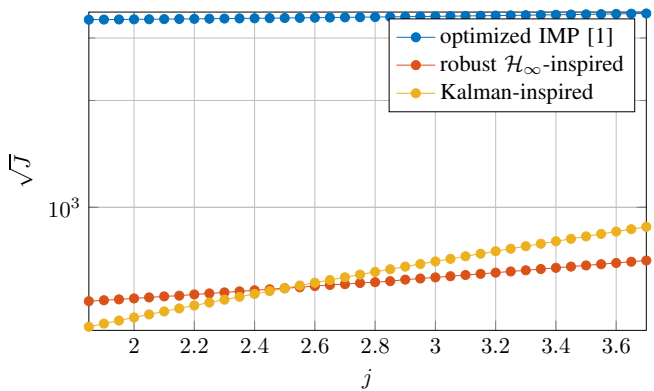
\begin{figure}[!ht]
\centering
\input{./Figures/figTest3.tex}
\caption{Square root of the cost $J$ from~\eqref{eq:Jcost} for an unstable $\F$, shown on a logarithmic scale as a function of measurement noise level $j$, comparing the performance of online gradient descent, robust $\mathcal{H}_{\infty}$-inspired, and Kalman-inspired algorithms.
}%
\label{fig:figTest3}%
\end{figure}%

In the first experiment, shown in Fig.~\ref{fig:figTest3}, we generate $\A$ as before, using the pair $\left(\F, \H\right)$ in an observable canonical form~\cite{joao_linear_2018}. Here, the characteristic polynomial of $\F$ is defined as $d_u\left(z\right)d_s\left(z\right)$, with $d_s\left(z\right)$ as in the experiment shown in Fig.~\ref{fig:figTest1} with $p=0.875$ and $d_u\left(z\right) = z^2 -2\cos\left(\omega_0\right)z + 1$, which is the $\mathcal{Z}$-transform of a sinusoidal signal with angular frequency $\omega_0 = \pi/12$. Using this setup, we repeat the experiment from Fig.~\ref{fig:figTest1}, varying $j$ in the interval $\left[1.85,\,3.7\right]$.
\begin{figure}[!ht]
\centering
\input{./Figures/figTest4.tex}
\caption{Square root of the cost $J$ from~\eqref{eq:Jcost} with an unstable $\F$, shown on a logarithmic scale as a function of $\lambda_{\max}$, comparing the performance of online gradient descent, robust $\mathcal{H}_{\infty}$-inspired, and Kalman-inspired algorithms.}%
\label{fig:figTest4}%
\end{figure}
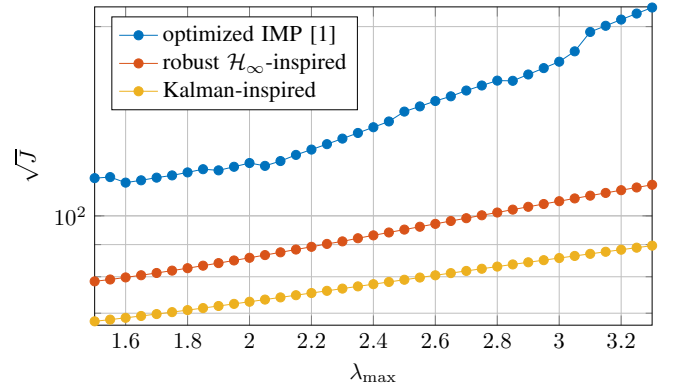%

The second experiment, depicted in Fig.~\ref{fig:figTest4}, replicates the test from Fig.~\ref{fig:figTest2} but for the case where $\F$ is unstable, as described in Fig.~\ref{fig:figTest3}, and with $j=1$. In this experiment, we vary $\lambda_{\max}$ within the interval $\left[1.5,\,3.3\right]$.

The simulation results presented in Figs.~\ref{fig:figTest3} and~\ref{fig:figTest4} illustrate that the proposed methods consistently outperform the optimized IMP algorithm from~\cite{bastianello_internal_2022}, both across different levels of measurement noise $j$ and as $\lambda_{\max}$ increases.

\section{Conclusions}\label{sec:conclusions}
In this paper, we introduced novel approaches for handling quadratic online optimization problems in which the minimum follows a noisy linear model and the quadratic cost is only partially known. This problem was addressed using two distinct control strategies: a Kalman-inspired algorithm and a robust $\mathcal{H}_\infty$-inspired algorithm. Both approaches demonstrated enhanced robustness and performance compared to conventional online gradient descent and other \textit{structured} algorithms. Future work may extend these methods to broader classes of online optimization problems where the quadratic assumption on the cost function is relaxed, and where other types of uncertainties can be addressed.
\bibliographystyle{IEEEtran}
\bibliography{references}
\end{document}

%% file: Figures/figTest0.tex
\resizebox{\columnwidth}{!}{%
\begin{tikzpicture}

\begin{axis}[%
scale only axis,
width = \columnwidth,
height = \columnwidth/1.75,
xmin=0,
xmax=5000,
xtick = {0,2500,5000},
minor x tick num=4,
xlabel={$k$},
xminorticks = true,
ymode=log,
ymin=10,
ymax=100,
ytick = {10,100},
yminorticks=true,
ylabel={$\norm{\e_k}$},
axis background/.style={fill=white},
xmajorgrids,
ymajorgrids,
yminorgrids,
legend style={legend cell align=left, align=left}
]
\addplot [color=mycolor1, line width=1.0pt]
  table[row sep=crcr]{%
0	45.2856105440664\\
100	49.4490766867856\\
200	54.1500430954856\\
300	53.4386044768606\\
400	54.7977635703322\\
500	55.5543739234667\\
600	54.6030620122799\\
700	56.2856726866335\\
800	56.8045755709844\\
900	56.1885406706279\\
1000	56.820964181078\\
1100	53.4909483538335\\
1200	52.7561404326784\\
1300	51.6558454102351\\
1400	51.3948717405249\\
1500	50.0190672704236\\
1600	50.983761726536\\
1700	50.7448965126795\\
1800	50.3900121240128\\
1900	52.4257264035497\\
2000	50.5017220124556\\
2100	50.8261800433686\\
2200	47.4534577217195\\
2300	53.449291509224\\
2400	54.0354473547992\\
2500	54.3665562313417\\
2600	54.0106578579622\\
2700	52.7166271746226\\
2800	52.2701154296169\\
2900	53.9657073312284\\
3000	55.7659383345325\\
3100	57.9290397573935\\
3200	60.0888106555726\\
3300	55.2734537562544\\
3400	53.2481728339063\\
3500	51.0005513277148\\
3600	49.1973532204672\\
3700	50.2018642791431\\
3800	52.3406802212955\\
3900	52.0798088123242\\
4000	53.8495410229653\\
4100	56.012851825609\\
4200	54.3466492520419\\
4300	55.0362628112578\\
4400	54.8182550798373\\
4500	56.623570155578\\
4600	58.5945512281178\\
4700	57.0222141629265\\
4800	55.5938401688194\\
4900	53.5708442026567\\
5000	49.937746014386\\
5100	46.7666548230235\\
5200	45.1973536067363\\
5300	43.6324909083485\\
5400	45.2066858108177\\
5500	43.1168267309085\\
5600	42.8830440349785\\
5700	44.6274802772551\\
5800	43.5630803278391\\
5900	40.7704504634706\\
6000	42.0352907519088\\
6100	41.1152353446723\\
6200	43.02005405146\\
6300	43.8989318348738\\
6400	43.742498746277\\
6500	47.9089406925678\\
6600	46.9328617077283\\
6700	46.8456074438431\\
6800	50.5105059483488\\
6900	52.5044887904982\\
7000	55.9610539393422\\
7100	55.4421111275563\\
7200	54.0264683190592\\
7300	55.4446652188155\\
7400	55.7969301418198\\
7500	54.8376351159571\\
7600	56.191961695124\\
7700	55.4097138661894\\
7800	51.8551433887624\\
7900	50.7935937869199\\
8000	46.273444285906\\
8100	46.7629315667029\\
8200	46.9088528629699\\
8300	44.1868616711268\\
8400	41.1326386131111\\
8500	41.656397768177\\
8600	42.0233594162695\\
8700	44.123202471173\\
8800	44.8183856362499\\
8900	44.7373730585711\\
9000	47.4377131685904\\
9100	48.9357367779362\\
9200	54.2352165650394\\
9300	57.0925230627682\\
9400	59.9594934500192\\
9500	57.7103298071672\\
9600	57.7560428999995\\
9700	57.6125671632562\\
9800	59.5490459879099\\
9900	62.6954874712368\\
};
\addlegendentry{gradient descent}

\addplot [color=mycolor2, line width=1.0pt]
  table[row sep=crcr]{%
0	16.7348196290301\\
100	16.9501185857881\\
200	16.5078826481993\\
300	16.3973456340494\\
400	16.4350963025194\\
500	16.342859931622\\
600	16.6142185403615\\
700	16.7369950237879\\
800	16.5941010425069\\
900	16.4537343680338\\
1000	16.587570742422\\
1100	16.4957519541606\\
1200	17.0952480095708\\
1300	16.8099765573298\\
1400	16.6186149819827\\
1500	16.7209147242776\\
1600	16.6820836615057\\
1700	16.5863172975266\\
1800	15.9519504432882\\
1900	15.9839540841365\\
2000	15.6308494140677\\
2100	15.5102205618428\\
2200	15.1268570353461\\
2300	15.4373936175889\\
2400	15.637760238843\\
2500	15.7125318230074\\
2600	15.6692925638408\\
2700	15.6911469956862\\
2800	15.7063025132532\\
2900	16.1103264640552\\
3000	16.7661703126141\\
3100	16.6332444550642\\
3200	16.3172159465\\
3300	16.5514259040582\\
3400	16.4359277157231\\
3500	16.3998256218873\\
3600	16.0972720002575\\
3700	16.2953880319092\\
3800	16.2441622538379\\
3900	16.0761931016373\\
4000	16.0458785640559\\
4100	16.4971983314131\\
4200	16.7348116365638\\
4300	16.6536328592107\\
4400	16.4341192206192\\
4500	16.1190532194748\\
4600	16.2281105880037\\
4700	16.7366041920619\\
4800	17.3761241627584\\
4900	17.8168215826572\\
5000	18.0766928955822\\
5100	17.8407492867226\\
5200	17.9853709262711\\
5300	17.6389477047249\\
5400	18.2883261128253\\
5500	18.4634707127677\\
5600	18.5420998525339\\
5700	17.7161461923695\\
5800	17.3809882579228\\
5900	16.8952169858389\\
6000	16.0210931307559\\
6100	16.1462050494898\\
6200	16.3288518322879\\
6300	16.5386248408735\\
6400	16.2054271794213\\
6500	16.2128653455518\\
6600	16.2947808823387\\
6700	16.0313764799278\\
6800	16.3075630275351\\
6900	16.1591072097358\\
7000	16.5575126339147\\
7100	16.2667793332555\\
7200	15.9046888271552\\
7300	16.0174714726332\\
7400	15.9104583159195\\
7500	16.0821017967685\\
7600	15.8836536796449\\
7700	16.4174864527405\\
7800	16.4054272189858\\
7900	16.4974711436024\\
8000	16.4544358708088\\
8100	16.6495457573197\\
8200	16.7510878404554\\
8300	16.3036645371605\\
8400	16.0347834385918\\
8500	16.0556005382409\\
8600	16.2479856964353\\
8700	16.060907829881\\
8800	15.588493299901\\
8900	15.365303470734\\
9000	15.3903656463125\\
9100	15.1755662828927\\
9200	15.3372154297662\\
9300	15.565603958344\\
9400	16.2847877457139\\
9500	16.086463655392\\
9600	16.0252554163371\\
9700	16.1048715908235\\
9800	16.1768137797887\\
9900	16.4280041292265\\
};
\addlegendentry{robust $\mathcal{H}_\infty$-inspired}

\addplot [color=mycolor3, line width=1.0pt]
  table[row sep=crcr]{%
0	17.402402403788\\
100	17.3186288225414\\
200	16.6737477844285\\
300	16.4745972887406\\
400	16.5487600762214\\
500	16.5556293632087\\
600	16.7558337253581\\
700	16.8145419389345\\
800	16.4591113053113\\
900	16.3705827448905\\
1000	16.4985103518126\\
1100	16.4030884566412\\
1200	16.975522265082\\
1300	16.7061664421768\\
1400	16.5150466491089\\
1500	16.5400111119399\\
1600	16.3739198775989\\
1700	16.2377132217642\\
1800	15.5790554322009\\
1900	15.5708652530432\\
2000	15.1448452655957\\
2100	15.0376686541115\\
2200	14.7755080298644\\
2300	15.0557258228285\\
2400	15.218100539624\\
2500	15.2731408263034\\
2600	15.3928291831936\\
2700	15.4815470928026\\
2800	15.3814516607667\\
2900	15.7374042976281\\
3000	16.4133719062723\\
3100	16.4481320811685\\
3200	16.120817232333\\
3300	16.5261833160697\\
3400	16.4719148719448\\
3500	16.6013384380753\\
3600	16.3674211167152\\
3700	16.6454710352571\\
3800	16.6415840732991\\
3900	16.4566800869139\\
4000	16.3441148189463\\
4100	16.6184102214508\\
4200	17.0405035678424\\
4300	16.9895221111749\\
4400	16.7483221161403\\
4500	16.2766515808737\\
4600	16.2265224976486\\
4700	16.6735935692737\\
4800	17.3959987851433\\
4900	17.9235969788815\\
5000	18.2757959642753\\
5100	18.1960017874304\\
5200	18.161925027051\\
5300	17.6213525359806\\
5400	18.3076897086693\\
5500	18.5462833754579\\
5600	18.6918208541295\\
5700	17.9347488623942\\
5800	17.6617728192297\\
5900	17.2043382535523\\
6000	16.3017199970897\\
6100	16.2677607452717\\
6200	16.4148098484382\\
6300	16.7263173451851\\
6400	16.2618964143955\\
6500	16.1203891826917\\
6600	16.2759340235216\\
6700	16.0296133262713\\
6800	16.0738584730051\\
6900	15.8433878372009\\
7000	16.2182056495066\\
7100	15.9793112469087\\
7200	15.8870823251101\\
7300	15.8834540473373\\
7400	15.7935925751981\\
7500	15.9827388853927\\
7600	15.8078080489835\\
7700	16.3527381098636\\
7800	16.5240955327678\\
7900	16.6393889128657\\
8000	16.8101635691938\\
8100	17.1988345708939\\
8200	17.0840291189837\\
8300	16.784504454429\\
8400	16.5754856180284\\
8500	16.4984729310276\\
8600	16.4768173663179\\
8700	16.2948018834656\\
8800	15.8370204620854\\
8900	15.6942760290404\\
9000	15.5663799001523\\
9100	15.272463256387\\
9200	15.3589259821739\\
9300	15.6347840465044\\
9400	16.3864859366015\\
9500	16.2312141544526\\
9600	16.3239184685322\\
9700	16.3244832754925\\
9800	16.3309031747683\\
9900	16.5124648029063\\
};
\addlegendentry{Kalman-inspired}

\end{axis}
\end{tikzpicture}%
}%

%% file: Figures/figTest1.tex
\resizebox{\columnwidth}{!}{%
\begin{tikzpicture}

\begin{axis}[%
scale only axis,
width = \columnwidth,
height = \columnwidth/1.75,
xmin=0.2,
xmax=2,
xlabel={$j$},
ymode=log,
ymin=100,
ymax=2125.89860948557,
ytick = {100,1000},
yminorticks=true,
ylabel={$\sqrt{J}$},
axis background/.style={fill=white},
xmajorgrids,
ymajorgrids,
yminorgrids,
legend style={legend cell align=left, align=left}
]
\addplot [color=mycolor1, mark=*, mark options={solid, fill=mycolor1, mycolor1}]
  table[row sep=crcr]{%
0.2	2105.7293996973\\
0.25	2104.52098501212\\
0.3	2103.41363721338\\
0.35	2102.40735630106\\
0.4	2101.50214227519\\
0.45	2100.69799513576\\
0.5	2099.99491488275\\
0.55	2099.3929015162\\
0.6	2098.89195503606\\
0.65	2098.49207544237\\
0.7	2098.19326273511\\
0.75	2097.99551691429\\
0.8	2097.89883797992\\
0.85	2097.90322593197\\
0.9	2098.00868077046\\
0.95	2098.21520249538\\
1	2098.52279110674\\
1.05	2098.93144660454\\
1.1	2099.44116898877\\
1.15	2100.05195825945\\
1.2	2100.76381441655\\
1.25	2101.57673746009\\
1.3	2102.49072739007\\
1.35	2103.50578420648\\
1.4	2104.62190790934\\
1.45	2105.83909849862\\
1.5	2107.15735597434\\
1.55	2108.5766803365\\
1.6	2110.0970715851\\
1.65	2111.71852972012\\
1.7	2113.4410547416\\
1.75	2115.26464664949\\
1.8	2117.18930544384\\
1.85	2119.21503112461\\
1.9	2121.34182369183\\
1.95	2123.56968314548\\
2	2125.89860948557\\
};
\addlegendentry{gradient descent}

\addplot [color=mycolor2,  mark=*, mark options={solid, fill=mycolor2, mycolor2}]
  table[row sep=crcr]{%
0.2	155.743085677622\\
0.25	157.450393317609\\
0.3	159.781129194453\\
0.35	161.574125149592\\
0.4	163.704867182762\\
0.45	166.761703123882\\
0.5	170.239866960935\\
0.55	173.322632791615\\
0.6	177.677668374798\\
0.65	181.832973665204\\
0.7	187.60707588351\\
0.75	192.874675142565\\
0.8	197.93589386104\\
0.85	204.085434036847\\
0.9	209.808925961692\\
0.95	216.265595989743\\
1	222.382790790428\\
1.05	229.001556248282\\
1.1	235.60668687286\\
1.15	242.506467910886\\
1.2	249.388278608314\\
1.25	256.7276087158\\
1.3	264.209022565776\\
1.35	271.358521842891\\
1.4	279.260775150247\\
1.45	286.829556302592\\
1.5	294.903630502801\\
1.55	302.626785591566\\
1.6	310.860116287594\\
1.65	319.306209068063\\
1.7	327.802682777126\\
1.75	336.242395442814\\
1.8	344.909949206595\\
1.85	352.671694924354\\
1.9	362.166700085264\\
1.95	371.734797425889\\
2	380.195727675549\\
};
\addlegendentry{robust $\mathcal{H}_\infty$-inspired}

\addplot [color=mycolor3,  mark=*, mark options={solid, fill=mycolor3, mycolor3}]
  table[row sep=crcr]{%
0.2	113.948977391187\\
0.25	118.935326683412\\
0.3	125.445292160552\\
0.35	132.335950738444\\
0.4	139.468469513992\\
0.45	146.800972449437\\
0.5	154.314752268572\\
0.55	161.999166476518\\
0.6	169.847094647475\\
0.65	177.853218646163\\
0.7	186.013259776565\\
0.75	194.323598002697\\
0.8	202.781063463541\\
0.85	211.382813431941\\
0.9	220.126255030522\\
0.95	229.008994061711\\
1	238.028799564081\\
1.05	247.183578285392\\
1.1	256.471355663452\\
1.15	265.890261229569\\
1.2	275.438517111695\\
1.25	285.114428770526\\
1.3	294.916377384233\\
1.35	304.842813477665\\
1.4	314.892251509916\\
1.45	325.063265213587\\
1.5	335.354483533336\\
1.55	345.764587049731\\
1.6	356.292304801348\\
1.65	366.936411438055\\
1.7	377.695724652815\\
1.75	388.569102850237\\
1.8	399.555443018435\\
1.85	410.653678776908\\
1.9	421.862778578179\\
1.95	433.181744044711\\
2	444.609608425628\\
};
\addlegendentry{Kalman-inspired}

\end{axis}
\end{tikzpicture}%
}%

%% file: Figures/figTest2.tex
\resizebox{\columnwidth}{!}{%
\begin{tikzpicture}

\begin{axis}[%
scale only axis,
width = \columnwidth,
height = \columnwidth/1.75,
xmin=2.55,
xmax=4.4,
xlabel={$\lambda_{\max}$},
ymode=log,
ymin=20,
ymax=700,
ytick={10,100,1000},
yminorticks=true,
ylabel={$\sqrt{J}$},
axis background/.style={fill=white},
xmajorgrids,
ymajorgrids,
yminorgrids,
legend style={legend cell align=left, align=left}
]
\addplot [color=mycolor1, mark=*, mark options={solid, fill=mycolor1, mycolor1}]
  table[row sep=crcr]{%
2.55	526.15782699729\\
2.6	530.236799845194\\
2.65	534.306183899321\\
2.7	538.36558395741\\
2.75	542.414643406899\\
2.8	546.453040658217\\
2.85	550.480485942424\\
2.9	554.496718431666\\
2.95	558.501503646214\\
3	562.494631116485\\
3.05	566.475912271683\\
3.1	570.445178531275\\
3.15	574.402279577125\\
3.2	578.347081787814\\
3.25	582.279466817853\\
3.3	586.199330307351\\
3.35	590.106580708673\\
3.4	594.001138218452\\
3.45	597.882933804654\\
3.5	601.751908319271\\
3.55	605.608011688204\\
3.6	609.451202171623\\
3.65	613.281445687179\\
3.7	617.098715190808\\
3.75	620.90299011003\\
3.8	624.694255823889\\
3.85	628.472503186528\\
3.9	632.237728089661\\
3.95	635.989931060617\\
4	639.729116893198\\
4.05	643.455294307946\\
4.1	647.168475639509\\
4.15	650.868676549105\\
4.2	654.555915759189\\
4.25	658.230214809047\\
4.3	661.891597829556\\
4.35	665.540091334994\\
4.4	669.17572403132\\
};
\addlegendentry{gradient descent}

\addplot [color=mycolor2, mark=*, mark options={solid, fill=mycolor2, mycolor2}]
  table[row sep=crcr]{%
2.55	23.3718390756647\\
2.6	23.7060153453466\\
2.65	24.041087623315\\
2.7	24.3531688698395\\
2.75	24.7029561204896\\
2.8	25.0413884615519\\
2.85	25.3743711752432\\
2.9	25.705831624259\\
2.95	26.0358932643973\\
3	26.37507614476\\
3.05	26.6966065455204\\
3.1	27.0240977455684\\
3.15	27.34697696659\\
3.2	27.6750976330549\\
3.25	28.0037149935637\\
3.3	28.3247125318845\\
3.35	28.6575325521352\\
3.4	28.9796702949453\\
3.45	29.3020592073543\\
3.5	29.6177621765236\\
3.55	29.9450391379644\\
3.6	30.2628030893161\\
3.65	30.5822059215474\\
3.7	30.8963557308442\\
3.75	31.2120760194735\\
3.8	31.5348936855931\\
3.85	31.843430958311\\
3.9	32.1563718432418\\
3.95	32.4689113896286\\
4	32.7908713815022\\
4.05	33.0921013563113\\
4.1	33.4084129825357\\
4.15	33.7213158282204\\
4.2	34.008876527305\\
4.25	34.333055643958\\
4.3	34.6373773158778\\
4.35	34.9429711405404\\
4.4	35.2492866593542\\
};
\addlegendentry{robust $\mathcal{H}_\infty$-inspired}

\addplot [color=mycolor3, mark=*, mark options={solid, fill=mycolor3, mycolor3}]
  table[row sep=crcr]{%
2.55	21.2944417887403\\
2.6	21.6656119212089\\
2.65	22.0456365022228\\
2.7	22.4347364772925\\
2.75	22.8331419713952\\
2.8	23.241090454027\\
2.85	23.6588247594939\\
2.9	24.0865909344392\\
2.95	24.5246358856115\\
3	24.9732048020735\\
3.05	25.4325383275501\\
3.1	25.9028694605552\\
3.15	26.3844201624474\\
3.2	26.8773976567902\\
3.25	27.3819904074805\\
3.3	27.8983637681919\\
3.35	28.42665530192\\
3.4	28.9669697769229\\
3.45	29.519373854239\\
3.5	30.0838904923174\\
3.55	30.6604931061409\\
3.6	31.2490995315336\\
3.65	31.8495658600234\\
3.7	32.4616802254825\\
3.75	33.0851566405078\\
3.8	33.7196289976636\\
3.85	34.3646453677919\\
3.9	35.0196627438628\\
3.95	35.6840423934425\\
4	36.3570459948984\\
4.05	37.0378327408098\\
4.1	37.7254575956623\\
4.15	38.4188708925863\\
4.2	39.1169194447201\\
4.25	39.8183493296929\\
4.3	40.5218104803001\\
4.35	41.2258631801957\\
4.4	41.9289865206462\\
};
\addlegendentry{Kalman-inspired}

\end{axis}
\end{tikzpicture}%
}%

%% file: Figures/figTest3.tex
\resizebox{\columnwidth}{!}{%
\begin{tikzpicture}

\begin{axis}[%
scale only axis,
width = \columnwidth,
height = \columnwidth/1.75,
xmin=1.85,
xmax=3.7,
xlabel={$j$},
ymode=log,
ymin=450,
ymax=3550,
ytick = {1000,10000},
yminorticks=true,
ylabel={$\sqrt{J}$},
axis background/.style={fill=white},
xmajorgrids,
ymajorgrids,
yminorgrids,
legend style={legend cell align=left, align=left}
]
\addplot [color=mycolor1,  mark=*, mark options={solid, fill=mycolor1, mycolor1}]
  table[row sep=crcr]{%
1.85	3378.67881346777\\
1.9	3381.10537154074\\
1.95	3383.60576635234\\
2	3386.17999790109\\
2.05	3388.82806618785\\
2.1	3391.54997121263\\
2.15	3394.34571297475\\
2.2	3397.21529147503\\
2.25	3400.15870671311\\
2.3	3403.17595868923\\
2.35	3406.26704740272\\
2.4	3409.43197285398\\
2.45	3412.67073504347\\
2.5	3415.98333397063\\
2.55	3419.36976963498\\
2.6	3422.83004203784\\
2.65	3426.364151178\\
2.7	3429.97209705683\\
2.75	3433.6538796726\\
2.8	3437.40949902601\\
2.85	3441.23895511843\\
2.9	3445.14224794711\\
2.95	3449.11937751522\\
3	3453.17034382007\\
3.05	3457.29514686304\\
3.1	3461.49378664391\\
3.15	3465.76626316212\\
3.2	3470.11257641845\\
3.25	3474.53272641287\\
3.3	3479.02671314464\\
3.35	3483.59453661449\\
3.4	3488.23619682198\\
3.45	3492.95169376802\\
3.5	3497.74102745036\\
3.55	3502.60419787128\\
3.6	3507.54120503053\\
3.65	3512.55204892703\\
3.7	3517.63672956126\\
};
\addlegendentry{optimized IMP~\cite{bastianello_internal_2022}}

\addplot [color=mycolor2,  mark=*, mark options={solid, fill=mycolor2, mycolor2}]
  table[row sep=crcr]{%
1.85	543.843820999201\\
1.9	547.121835589972\\
1.95	550.240302274121\\
2	553.635930147404\\
2.05	557.31902475688\\
2.1	561.011715482071\\
2.15	563.55792850982\\
2.2	568.186811543581\\
2.25	571.966342357822\\
2.3	575.72952551788\\
2.35	579.629339047065\\
2.4	583.693485175819\\
2.45	587.495013980424\\
2.5	591.554879029277\\
2.55	595.063173159574\\
2.6	599.957748659371\\
2.65	604.19445490691\\
2.7	607.368388258485\\
2.75	612.008241332869\\
2.8	615.206440323474\\
2.85	619.251122657469\\
2.9	625.453014225615\\
2.95	630.066061774428\\
3	634.637548388676\\
3.05	639.471478824444\\
3.1	644.477150042412\\
3.15	649.380458419065\\
3.2	654.403966242681\\
3.25	659.604687307009\\
3.3	664.539031251354\\
3.35	670.034197684565\\
3.4	675.249497786504\\
3.45	680.583228854137\\
3.5	685.863675477964\\
3.55	691.386944097696\\
3.6	696.907329862908\\
3.65	702.522841076772\\
3.7	708.149180054473\\
};
\addlegendentry{robust $\mathcal{H}_\infty$-inspired}

\addplot [color=mycolor3,  mark=*, mark options={solid, fill=mycolor3, mycolor3}]
  table[row sep=crcr]{%
1.85	460.714146114936\\
1.9	470.14267242373\\
1.95	479.672280755411\\
2	489.303335467776\\
2.05	499.036224735106\\
2.1	508.871359496757\\
2.15	518.809172558339\\
2.2	528.850117822869\\
2.25	538.994669636683\\
2.3	549.243322229986\\
2.35	559.596589244798\\
2.4	570.055003334642\\
2.45	580.619115831226\\
2.5	591.289496470391\\
2.55	602.066733166681\\
2.6	612.951431838016\\
2.65	623.944216268956\\
2.7	635.045728016201\\
2.75	646.256626345022\\
2.8	657.577588200265\\
2.85	669.00930820574\\
2.9	680.552498692694\\
2.95	692.207889751355\\
3	703.976229308579\\
3.05	715.85828322757\\
3.1	727.854835428662\\
3.15	739.966688034747\\
3.2	752.194661527591\\
3.25	764.539594932194\\
3.3	777.00234601507\\
3.35	789.583791502022\\
3.4	802.284827310613\\
3.45	815.106368804374\\
3.5	828.049351059698\\
3.55	841.11472915108\\
3.6	854.303478453089\\
3.65	867.616594957944\\
3.7	881.055095612596\\
};
\addlegendentry{Kalman-inspired}

\end{axis}
\end{tikzpicture}%
}%

%% file: Figures/figTest4.tex
\resizebox{\columnwidth}{!}{%
\begin{tikzpicture}
\begin{axis}[%
scale only axis,
width = \columnwidth,
height = \columnwidth/1.75,
xmin=1.5,
xmax=3.3,
xlabel={$\lambda_{\max}$},
ymode=log,
ymin=67,
ymax=215,
ytick={10,100, 1000},
yminorticks=true,
ylabel={$\sqrt{J}$},
axis background/.style={fill=white},
xmajorgrids,
ymajorgrids,
yminorgrids,
legend style={legend cell align=left, align=left},
legend pos = north west
]
\addplot [color=mycolor1, mark=*, mark options={solid, fill=mycolor1, mycolor1}]
  table[row sep=crcr]{%
1.5	114.81745306673\\
1.55	115.257690430443\\
1.6	112.942997000836\\
1.65	113.903461729995\\
1.7	115.022057962307\\
1.75	115.96590289284\\
1.8	117.250421002768\\
1.85	118.630226609935\\
1.9	118.13538181636\\
1.95	119.687797454699\\
2	121.34665475124\\
2.05	120.086169118182\\
2.1	122.225730072131\\
2.15	125.010448720127\\
2.2	127.467743978778\\
2.25	130.03249284317\\
2.3	132.702058238711\\
2.35	135.472962927005\\
2.4	138.347081987193\\
2.45	141.305139037888\\
2.5	146.507620512771\\
2.55	149.321963853154\\
2.6	152.306140450468\\
2.65	154.945577514238\\
2.7	158.257324049668\\
2.75	161.202579542115\\
2.8	164.113491149681\\
2.85	164.00022461207\\
2.9	167.844672479703\\
2.95	171.807387608568\\
3	175.883638000439\\
3.05	182.684867795917\\
3.1	196.150649696301\\
3.15	200.612130343625\\
3.2	205.225535692946\\
3.25	209.805326174595\\
3.3	214.487210184892\\
};
\addlegendentry{optimized IMP~\cite{bastianello_internal_2022}}

\addplot [color=mycolor2,  mark=*, mark options={solid, fill=mycolor2, mycolor2}]
  table[row sep=crcr]{%
1.5	78.6673399083772\\
1.55	79.2186015941112\\
1.6	79.8175113446943\\
1.65	80.4452138054556\\
1.7	81.1178031141282\\
1.75	81.818887403974\\
1.8	82.553908832468\\
1.85	83.3183099471598\\
1.9	84.1073719239918\\
1.95	84.9221736137991\\
2	85.7590881813219\\
2.05	86.6154816092443\\
2.1	87.4936793605255\\
2.15	88.3893657405033\\
2.2	89.2764053820723\\
2.25	90.2305159931397\\
2.3	91.0709138796491\\
2.35	92.1287477638569\\
2.4	93.1002474213276\\
2.45	94.0894535164685\\
2.5	95.0877445940048\\
2.55	96.0888241155553\\
2.6	97.1055765345825\\
2.65	98.1303546864439\\
2.7	99.1407387273859\\
2.75	100.202938022808\\
2.8	101.229807014753\\
2.85	102.306840183894\\
2.9	103.368897994363\\
2.95	104.432158745027\\
3	105.511173824218\\
3.05	106.590910860307\\
3.1	107.674175734345\\
3.15	108.7621567938\\
3.2	109.853450564804\\
3.25	110.947604883723\\
3.3	112.047758404565\\
};
\addlegendentry{robust $\mathcal{H}_\infty$-inspired}

\addplot [color=mycolor3,  mark=*, mark options={solid, fill=mycolor3, mycolor3}]
  table[row sep=crcr]{%
1.5	67.95936315185\\
1.55	68.3837083430002\\
1.6	68.8304341843399\\
1.65	69.2976204337884\\
1.7	69.7835253938672\\
1.75	70.2865655106508\\
1.8	70.8052977423559\\
1.85	71.3384042681113\\
1.9	71.8846791816818\\
1.95	72.4430168750172\\
2	73.0124018654043\\
2.05	73.5918998600301\\
2.1	74.1806498846191\\
2.15	74.7778573299228\\
2.2	75.3827877923086\\
2.25	75.9947616033563\\
2.3	76.6131489589662\\
2.35	77.2373655715203\\
2.4	77.8668687795976\\
2.45	78.5011540589763\\
2.5	79.1397518864789\\
2.55	79.7822249147841\\
2.6	80.4281654220048\\
2.65	81.0771930045779\\
2.7	81.7289524861181\\
2.75	82.3831120183844\\
2.8	83.0393613535244\\
2.85	83.697410269332\\
2.9	84.3569871315015\\
2.95	85.0178375787584\\
3	85.6797233184632\\
3.05	86.3424210216856\\
3.1	87.0057213080283\\
3.15	87.6694278115852\\
3.2	88.3333563203635\\
3.25	88.9973339823345\\
3.3	89.6611985720282\\
};
\addlegendentry{Kalman-inspired}

\end{axis}
\end{tikzpicture}%
}%

%% file: main-new.bbl
\begin{thebibliography}{10}
\providecommand{\url}[1]{#1}
\csname url@samestyle\endcsname
\providecommand{\newblock}{\relax}
\providecommand{\bibinfo}[2]{#2}
\providecommand{\BIBentrySTDinterwordspacing}{\spaceskip=0pt\relax}
\providecommand{\BIBentryALTinterwordstretchfactor}{4}
\providecommand{\BIBentryALTinterwordspacing}{\spaceskip=\fontdimen2\font plus
\BIBentryALTinterwordstretchfactor\fontdimen3\font minus
  \fontdimen4\font\relax}
\providecommand{\BIBforeignlanguage}[2]{{%
\expandafter\ifx\csname l@#1\endcsname\relax
\typeout{** WARNING: IEEEtran.bst: No hyphenation pattern has been}%
\typeout{** loaded for the language `#1'. Using the pattern for}%
\typeout{** the default language instead.}%
\else
\language=\csname l@#1\endcsname
\fi
#2}}
\providecommand{\BIBdecl}{\relax}
\BIBdecl

\bibitem{bastianello_internal_2022}
N.~Bastianello, R.~Carli, and S.~Zampieri, ``Internal {Model}-{Based} {Online}
  {Optimization},'' \emph{IEEE Transactions on Automatic Control}, pp. 1--8,
  2023.

\bibitem{liao-mcpherson_semismooth_2018}
D.~Liao-McPherson, M.~Nicotra, and I.~Kolmanovsky, ``A {Semismooth} {Predictor}
  {Corrector} {Method} for {Real}-{Time} {Constrained} {Parametric}
  {Optimization} with {Applications} in {Model} {Predictive} {Control},'' in
  \emph{2018 {IEEE} {Conference} on {Decision} and {Control} ({CDC})}, Dec.
  2018, pp. 3600--3607.

\bibitem{paternain_realtime_2019}
S.~Paternain, M.~Morari, and A.~Ribeiro, ``Real-{Time} {Model} {Predictive}
  {Control} {Based} on {Prediction}-{Correction} {Algorithms},'' in \emph{2019
  {IEEE} 58th {Conference} on {Decision} and {Control} ({CDC})}, 2019, pp.
  5285--5291.

\bibitem{hall_online_2015}
E.~C. Hall and R.~M. Willett, ``Online {Convex} {Optimization} in {Dynamic}
  {Environments},'' \emph{IEEE Journal of Selected Topics in Signal
  Processing}, vol.~9, no.~4, pp. 647--662, Jun. 2015.

\bibitem{fosson_centralized_2021}
S.~M. Fosson, ``Centralized and {Distributed} {Online} {Learning} for {Sparse}
  {Time}-{Varying} {Optimization},'' \emph{IEEE Transactions on Automatic
  Control}, vol.~66, no.~6, pp. 2542--2557, Jun. 2021.

\bibitem{natali_online_2021}
A.~Natali, M.~Coutino, E.~Isufi, and G.~Leus, ``Online {Time}-{Varying}
  {Topology} {Identification} {Via} {Prediction}-{Correction} {Algorithms},''
  in \emph{{ICASSP} 2021 - 2021 {IEEE} {International} {Conference} on
  {Acoustics}, {Speech} and {Signal} {Processing} ({ICASSP})}.\hskip 1em plus
  0.5em minus 0.4em\relax Toronto, ON, Canada: IEEE, Jun. 2021, pp. 5400--5404.

\bibitem{shalev_online_2011}
S.~Shalev-Shwartz, ``Online {Learning} and {Online} {Convex} {Optimization},''
  \emph{Foundations and Trends® in Machine Learning}, vol.~4, no.~2, pp.
  107--194, 2011.

\bibitem{dixit_online_2019}
R.~Dixit, A.~S. Bedi, R.~Tripathi, and K.~Rajawat, ``Online {Learning} with
  {Inexact} {Proximal} {Online} {Gradient} {Descent} {Algorithms},'' \emph{IEEE
  Transactions on Signal Processing}, vol.~67, no.~5, pp. 1338 -- 1352, 2019.

\bibitem{chang_distributed_2020}
T.-H. Chang, M.~Hong, H.-T. Wai, X.~Zhang, and S.~Lu, ``Distributed {Learning}
  in the {Nonconvex} {World}: {From} batch data to streaming and beyond,''
  \emph{IEEE Signal Processing Magazine}, vol.~37, no.~3, pp. 26--38, May 2020.

\bibitem{dallanese_optimization_2020}
E.~Dall'Anese, A.~Simonetto, S.~Becker, and L.~Madden, ``Optimization and
  {Learning} {With} {Information} {Streams}: {Time}-varying algorithms and
  applications,'' \emph{IEEE Signal Processing Magazine}, vol.~37, no.~3, pp.
  71--83, May 2020.

\bibitem{simonetto_timevarying_2020}
A.~Simonetto, E.~Dall'Anese, S.~Paternain, G.~Leus, and G.~B. Giannakis,
  ``Time-{Varying} {Convex} {Optimization}: {Time}-{Structured} {Algorithms}
  and {Applications},'' \emph{Proceedings of the IEEE}, vol. 108, no.~11, pp.
  2032--2048, Nov. 2020.

\bibitem{simonetto_dual_2019}
A.~Simonetto, ``Dual {Prediction}–{Correction} {Methods} for {Linearly}
  {Constrained} {Time}-{Varying} {Convex} {Programs},'' \emph{IEEE Transactions
  on Automatic Control}, vol.~64, no.~8, pp. 3355--3361, Aug. 2019.

\bibitem{bastianello_extrapolation_2023}
N.~Bastianello, R.~Carli, and A.~Simonetto, ``Extrapolation-{Based}
  {Prediction}-{Correction} {Methods} for {Time}-varying {Convex}
  {Optimization},'' \emph{Signal Processing}, vol. 210, p. 109089, 2023.

\bibitem{fadali_digital_2019}
M.~S. Fadali and A.~Visioli, \emph{Digital control engineering: analysis and
  design}, 3rd~ed.\hskip 1em plus 0.5em minus 0.4em\relax San Diego: Academic
  press is an imprint of Elsevier, 2019.

\bibitem{lessard_analysis_2016}
L.~Lessard, B.~Recht, and A.~Packard, ``Analysis and {Design} of {Optimization}
  {Algorithms} via {Integral} {Quadratic} {Constraints},'' \emph{SIAM Journal
  on Optimization}, vol.~26, no.~1, pp. 57--95, Jan. 2016.

\bibitem{scherer_optimization_2023}
C.~W. Scherer, C.~Ebenbauer, and T.~Holicki, ``Optimization {Algorithm}
  {Synthesis} based on {Integral} {Quadratic} {Constraints}: {A} {Tutorial},''
  Jun. 2023.

\bibitem{davydov_contracting_2023}
\BIBentryALTinterwordspacing
A.~Davydov, V.~Centorrino, A.~Gokhale, G.~Russo, and F.~Bullo, ``Contracting
  {Dynamics} for {Time}-{Varying} {Convex} {Optimization},'' May 2023,
  arXiv:2305.15595 [cs, eess, math]. [Online]. Available:
  \url{http://arxiv.org/abs/2305.15595}
\BIBentrySTDinterwordspacing

\bibitem{van_control_2021}
J.~H. van Schuppen, \emph{Control and system theory of discrete-time stochastic
  systems}.\hskip 1em plus 0.5em minus 0.4em\relax Springer, 2021.

\bibitem{debruyne_linear_1995}
\BIBentryALTinterwordspacing
F.~{De Bruyne}, B.~Anderson, and M.~Gevers, ``Relating $\mathcal{H}_2$ and
  $\mathcal{H}_\infty$ bounds for finite-dimensional systems,'' \emph{Systems
  \& Control Letters}, vol.~24, no.~3, pp. 173--181, 1995. [Online]. Available:
  \url{https://www.sciencedirect.com/science/article/pii/016769119400018Q}
\BIBentrySTDinterwordspacing

\bibitem{zhou_robust_1998}
K.~Zhou and J.~C. Doyle, \emph{Essentials of Robust Control}.\hskip 1em plus
  0.5em minus 0.4em\relax Prentice-Hall, 1998.

\bibitem{scherer_theory_2001}
C.~Scherer, ``Theory of robust control,'' \emph{Delft University of
  Technology}, pp. 1--160, 2001.

\bibitem{joao_linear_2018}
\BIBentryALTinterwordspacing
J.~P. Hespanha, \emph{Linear Systems Theory: Second Edition}, ned - new
  edition, 2~ed.\hskip 1em plus 0.5em minus 0.4em\relax Princeton University
  Press, 2018. [Online]. Available:
  \url{http://www.jstor.org/stable/j.ctvc772kp}
\BIBentrySTDinterwordspacing

\end{thebibliography}
